\documentclass[a4paper]{article}
\usepackage{amsthm,amsmath,amssymb,color}
\usepackage{color}

\newtheorem{teor}{Theorem}[section]
\newtheorem{defin}[teor]{Definition}
\newtheorem{lemm}[teor]{Lemma}
\newtheorem{osse}[teor]{Remark}
\newtheorem{prop}[teor]{Proposition}
\newtheorem{defi}[teor]{Definition}
\newtheorem{coro}[teor]{Corollary}
\newtheorem{prob}[teor]{Problem}

\newcommand{\bele}{\begin{lemm}\begin{sl}}
\newcommand{\enle}{\end{sl}\end{lemm}}
\newcommand{\bedef}{\begin{defi}\begin{sl}}
\newcommand{\eddef}{\end{sl}\end{defi}}
\newcommand{\bete}{\begin{teor}\begin{sl}}
\newcommand{\ente}{\end{sl}\end{teor}}
\newcommand{\beos}{\begin{osse}\begin{rm}}
\newcommand{\eddos}{\end{rm}\end{osse}}
\newcommand{\bepr}{\begin{prop}\begin{sl}}
\newcommand{\empr}{\end{sl}\end{prop}}
\newcommand{\bepro}{\begin{prob}\begin{rm}}
\newcommand{\empro}{\end{rm}\end{prob}}
\newcommand{\bede}{\begin{defin}\begin{sl}}
\newcommand{\edde}{\end{sl}\end{defin}}
\newcommand{\beco}{\begin{coro}\begin{sl}}
\newcommand{\enco}{\end{sl}\end{coro}}
\newcommand{\disp}{\displaystyle}

\newcommand{\qquand}{\qquad\text{and}\qquad}
\newcommand{\quext}{\quad\text}
\newcommand{\qquext}{\qquad\text}
\newcommand{\de}{\partial}
\newcommand{\Div}{\hbox{div}}

\newcommand{\RR}{\mathbb{R}}

\newcommand{\EE}{\mathbb{E}}

\newcommand{\NN}{\mathbb{N}}

\newcommand{\beeq}[1]{\begin{equation}\label{#1}}
\newcommand{\eddeq}{\end{equation}}

\newcommand{\beeqa}[1]{\begin{eqnarray}\label{#1}}
\newcommand{\eddeqa}{\end{eqnarray}}

\newcommand{\beal}[1]{\begin{align}\label{#1}}
\newcommand{\eddal}{\end{align}}

\newcommand{\bespl}[1]{\begin{split}\label{#1}}
\newcommand{\edspl}{\end{split}}

\newcommand{\bega}[1]{\begin{gather}\label{#1}}
\newcommand{\edga}{\end{gather}}

\newcommand{\beeqax}{\begin{eqnarray*}}
\newcommand{\eddeqax}{\end{eqnarray*}}

\def\qed{\ifmmode % if math mode, assume display: omit penalty etc.
  \else \leavevmode\unskip\penalty9999 \hbox{}\nobreak\hfill
  \fi
  \quad\hbox{\hskip.5em\vrule width.4em height.6em depth.05em\hskip.1em}}
\def\endproofsym{\qed}
\renewenvironment{proof}[1][Proof]{\trivlist\item[\hskip\labelsep{\hskip0pt
    {\normalfont\scshape#1.}\hskip .321429\parindent}]\ignorespaces}
{\endproofsym\endtrivlist}
\def\endnobox{\def\endproofsym{}\end{proof}\def\endproofsym{\qed}}

\newcommand{\no}{\nonumber}

\newcommand{\beeqao}{\begin{eqnarray}\no}
\newcommand{\bealo}{\begin{align}\no}
\newcommand{\besplo}{\begin{split}\no}
\newcommand{\begao}{\begin{gather}\no}

\newcommand{\ttt}{{\mathfrak t}}

\newcommand{\perogni}{\forall\,}
\newcommand{\esiste}{\exists\,}
\newcommand{\itt}{\int_0^t}

\newcommand{\io}{\int_\Omega}

\newcommand{\iga}{\int_\Gamma}

\newcommand{\ito}{\itt\!\io}
\newcommand{\itga}{\itt\!\iga}

\newcommand{\epsi}{\varepsilon}
\newcommand{\ee}{_{\varepsilon}}
\newcommand{\Ga}{_{\Gamma}}

\newcommand{\OO}{_{\Omega}}

\newcommand{\bn}{\boldsymbol{n}}

\newcommand{\dn}{\partial_{\bn}}
\newcommand{\fhi}{\varphi}

\newcommand{\lhs}{left hand side}
\newcommand{\rhs}{right hand side}

\DeclareMathOperator{\dive}{div}
\DeclareMathOperator{\deriv}{d}

\DeclareMathOperator{\sign}{sign}

\DeclareMathOperator{\loc}{loc}

\let\TeXchi\chi
\def\chi{{\setbox0 \hbox{\mathsurround0pt
$\TeXchi$}\hbox{\raise\dp0 \copy0 }}}

\newcommand{\znn}{_{0,n}}

\newcommand{\zzn}{_{0,n}}

\newcommand{\zzm}{_{0,m}}
\newcommand{\teta}{\vartheta}

\newcommand{\calH}{{\mathcal H}}

\newcommand{\calT}{{\mathcal T}}

\newcommand{\calE}{{\mathcal E}}

\newcommand{\calR}{{\mathcal R}}

\newcommand{\calV}{{\mathcal V}}

\newcommand{\calI}{{\mathcal I}}

\newcommand{\barO}{\overline{\Omega}}

\newcommand{\tetasu}{\overline{\teta}}
\newcommand{\tetagiu}{\underline{\teta}}

\newcommand{\dit}{\deriv\!t}

\newcommand{\dix}{\deriv\!x}

\newcommand{\dm}{\deriv\!m}
\newcommand{\ibaro}{\int_{\overline{\Omega}}}

\newcommand{\ddt}{\frac{\deriv\!{}}{\dit}}

\definecolor{coloras}{rgb}{0.,0.67,0}
\newenvironment{giuliorev}{\color{red}}{\color{black}}
\newcommand{\III}{\begin{giuliorev}}
\newcommand{\EEE}{\end{giuliorev}}

\numberwithin{equation}{section}
\begin{document}

\title{On a singular heat equation with dynamic boundary
  conditions}

\author{Giulio Schimperna\\
Dipartimento di Matematica, Universit\`a di Pavia,\\
Via Ferrata~1, I-27100 Pavia, Italy\\
E-mail: {\tt giusch04@unipv.it}\\
\and
Antonio Segatti\\
Dipartimento di Matematica, Universit\`a di Pavia,\\
Via Ferrata~1, I-27100 Pavia, Italy\\
E-mail: {\tt antonio.segatti@unipv.it}\\
\and
Sergey Zelik\\
Department of Mathematics, University of Surrey,\\
Guildford, GU2 7XH, United Kingdom\\
E-mail: {\tt S.Zelik@surrey.ac.uk}
}

%\date{}
\maketitle
\begin{abstract}
 In this paper we analyze a singular heat equation
 of the form $\teta_t + \Delta \teta^{-1} = f$.
 The singular term $\teta^{-1}$ gives rise to
 very fast diffusion effects. The equation is
 settled in a smooth bounded domain
 $\Omega\subset\RR^3$ and complemented
 with a general {\sl dynamic}\/ boundary condition of the
 form $\alpha \teta_t - \beta \Delta\Ga \teta = \dn \teta^{-1}$,
 where $\Delta\Ga$ is the Laplace-Beltrami operator
 and $\alpha$ and $\beta$ are nonnegative coefficients
 (in particular, the homogeneous Neumann case given by
 $\alpha=\beta=0$ is included). For this
 problem, we first introduce a suitable weak formulation
 and prove a related existence result. For more regular initial
 data, we show that there exists at least one
 weak solution satisfying instantaneous regularization effects
 which are uniform with respect to the time variable.
 In this improved regularity class, uniqueness
 is also shown to hold.
\end{abstract}


\vspace{2mm}

\noindent {\bf AMS (MOS) subject clas\-si\-fi\-ca\-tion:}~~%
 35K55, 35K67, 35K51, 80A22.

%%%%%%%%%%%%%%%%%%%%%%%%%%%%%%%%%%%%%%%%%%%%%%%%%%%%%%%%%%%%%%%%%%%
%%%%%%%%%%%%%%%%%%%%%% Caption da ambo i lati %%%%%%%%%%%%%%%%%%%%%

%
%
%\pagestyle{myheadings}
%\newcommand\testopari{\sc Giulio~Schimperna}
%\newcommand\testodispari{\sc Cahn-Hilliard equation}
%\markboth{\testodispari}{\testopari}

%%%%%%%%%%%%%%%%%%%%%%%%%%%%%%%%%%%%%%%%%%%%%%%%%%%%%%%%%%%%%%%%%%%
%%%%%%%%%%%%%%%%%%%%%%%%%%%%%%%%%%%%%%%%%%%%%%%%%%%%%%%%%%%%%%%%%%%

\section{Introduction}

In this paper we address a singular heat equation
having the expression
\begin{equation}\label{calorein}
  \teta_t + \Delta \teta^{-1} = f,
\end{equation}
where $f$ is an external force. The equation
is settled in a smooth bounded domain $\Omega\subset\RR^3$,
where the restriction to the three-dimensional setting
is motivated by physical applications, and is complemented with
{\sl dynamic}\/ boundary conditions of the form
\begin{equation}\label{dyn:in}
  \alpha \teta_t - \beta \Delta\Ga \teta = \dn \teta^{-1}
   \quext{on }\,\Gamma:=\de\Omega,
\end{equation}
where $\Delta\Ga$ is the Laplace-Beltrami operator and
$\alpha,\beta\ge 0$ (in particular they may both
be $0$, so that the homogeneous Neumann problem is included).
%Condition \eqref{dyn:in} ensures conservation of the total
%``mass'' of $\teta$ (this includes also the boundary mass
%if $\alpha>0$).
Equation \eqref{calorein} can be viewed
in the framework of nonlinear diffusion problems, i.e.,
of equations of the form
\begin{equation}\label{eq:diffusion}
  \teta_t -\Div(\teta^{m-1}\nabla \teta) = f,
\end{equation}
with $m\in\mathbb{R}$. When $m\ge1$, equation \eqref{eq:diffusion} is the
well-known porous medium equation \cite{Va07} (heat equation when $m=1$);
for $m<1$, equation \eqref{eq:diffusion} lies in the class of fast diffusion
equations (ultra fast when $m<0$) analyzed in several papers 
(see, e.g., the pioneering work~\cite{HP}).
In this regime very fast diffusion occurs in the regions where $\teta$ 
is small. Moreover, one can distinguish two sub-regimes, depending on  
the so-called {\sl first critical fast diffusion exponent}\/ $m_c:=\frac{d-2}{2}$
($d$ is the dimension) which acts as a threshold between the good
parameter range $m\in(m_c,1)$ and its complementary range $m\le m_c$.
More precisely, when $m_c < m < 1$ and the initial condition
$\teta_0$ is non-negative and locally integrable, then equation
\eqref{eq:diffusion} (with $f=0$) on the whole space has a unique
weak solution which is globally defined, positive and smooth.
Moreover, solutions emanating from initial data lying in
$L^p(\mathbb{R}^d)$ ($p\ge 1$), or even in the Marcinkiewicz
space $M^p(\mathbb{R}^d)$, immediately become bounded. The
scenario when $m < m_c$ is drastically different
for at least two reasons (we refer to \cite{BonVaz09} for further remarks).
First of all, solutions may in general be unbounded and
may also be non-smooth. As an example, when
$\Omega=\mathbb{R}^3$, one can consider, for arbitrary $T>0$,
the function (see \cite{Va06})
\begin{equation} \label{explicit}
  \Theta(t,x)=2\big((T-t)^{+}\vert x\vert^{-2}\big)^{1/2},
\end{equation}
which solves \eqref{calorein} with $f=0$ and initial condition
$\Theta_0(x)=2 T^{1/2}\vert x\vert^{-1}\in L^{p}_{\loc}(\mathbb{R}^3)$
for any $p<3$. In particular, $\Theta$ is an unbounded function and
$\Theta(t,\cdot)\notin L^{p}_{\loc}(\mathbb{R}^3)$ for $p\ge 3$, at
least until it vanishes for $t\ge T$.
In this range of $m$, the boundedness of the solutions is tied to
the summability of the initial condition and to the value of $m$.
More precisely, solutions are bounded whenever
$\teta_0 \in L^p(\mathbb{R}^d)$ with $p > p_c:=\frac{d(1-m)}{2}$
(see \cite{BonVaz09}). Note that, when $d=3$ and $m=-1$
(as in our case \eqref{calorein}),
$p_c = 3$. A second notable feature of the fast diffusion regime,
which is already evident in example \eqref{explicit},
is the possible occurrence of extinction in finite time.
This means that a solution may become identically zero
after some finite time $T$ which depends on the initial conditions.
Consequently, positivity is lost. The ultra fast diffusion
regime $m < 0$ presents further difficulties linked
to the mere question of existence (see \cite{Vaz92}
and \cite{DeskaDelp97}). In particular, for
(nonzero) data in $L^p$ with $p <p_c$
one may face a phenomenon called ``immediate extinction''
meaning that solutions obtained as limits of reasonable
approximation schemes can be identically zero for any
$t > 0$. Note that the immediate extinction can occur also
for boundary value problems with zero Dirichlet conditions.
Finally, the case $m = m_* = \frac{d-4}{d-2}$, $d\ge 3$,
deserves a particular attention, as observed in
the papers \cite{BlaBonDolGrVaz09}
and \cite{BGV-arma10} dealing with the asymptotics
as $t\nearrow T$ (the extinction time). Note that we
always have $m_* < m_c$ and, for $d=3$, we have $m_* =-1$,
exactly as in our equation \eqref{calorein}.
%
%The literature on equations of the form
%\eqref{eq:diffusion} in the case $m>0$ is very large (see, e.g.,
%the book \cite{Va06}).
%%\cite{FK}, {(\bf here, quote some Vazquez or
%%Grillo or look for some survey paper)}, see also \cite{BG}).
%%In the case $m\in(0,1)$ (fast diffusion equation), \dots\dots
%On the other hand, the literature on the ultra-fast
% regime $m<0$ is far less wide. Most of the results are connected with
%the whole space situation.
%In \cite{Hui}, the author shows existence and uniqueness
%of solutions for the case $m\in (-1,0)$ with zero source $f$,
%but with nonhomogeneous Neumann conditions. {\bf Say something
%else}.
%\\

In this paper, we are interested in the analysis of equation
\eqref{calorein} with the dynamic boundary condition
\eqref{dyn:in} in the case when the source term $f$
has zero spatial mean. Then, a straightforward
computation permits us to see that conservation of mass occurs.
For instance, in the case $\alpha  = 0$, we have
$\int_\Omega\teta(t,x)\dix = \int_\Omega\teta_0(x)\dix$
for any $t>0$, which manifestly excludes both immediate
extinction and extinction in finite time. 
It is worth pointing out that Neumann (or related)
boundary problems are poorly studied in the literature, even
in the simpler case of the porous medium equation. In this regard,
well-posedness and asymptotic results have been established
in~\cite{AR}. More recently, these results have been improved
in~\cite{GM} by means of an approach similar to ours based on
Moser iteration techniques. 

Our interest is twofold. From the one hand, we
aim at proving existence of at least one solution
under weak conditions on the initial data.
To explain what we mean for ``weak'', and considering for
simplicity the case $f=0$, we note that, thanks to the
monotone structure of \eqref{calorein}-\eqref{dyn:in},
the system admits several Liapunov
functionals. In particular, on account of 
physical considerations (see below for details),
we may identify a natural {\sl energy}\/ functional $\calE$,
which is defined in \eqref{defiE} below. Mathematically,
the finiteness of $\calE$ seems to be a minimal
regularity condition on the initial datum $\teta_0$
that allows to define a rigorous concept of weak solution
(cf.~Definition~\ref{def:weak} below). Indeed
(in the simpler Neumann case $\alpha=\beta=0$),
finiteness of the energy corresponds to asking that
$\teta_0-\log\teta_0\in L^1(\Omega)$. Hence, we
have some control on the $L^1$-norm of the
solution (permitting to use $L^1$-arguments in the
analysis) and also a positivity condition.
Weak solutions emanating from initial data
with finite energy $\calE$ will be called
``energy'' solutions. Actually, our first result
(see Theorem \ref{teo:energy} below) states that,
if the initial data have finite energy
and $f$ has zero mean value and
satisfies suitable summability conditions,
then at least one global energy solution exists;
moreover, the energy $\calE(\teta)$ remains bounded
uniformly in~time.

Once existence is established, we study boundedness and
positivity properties of energy solutions. Assuming the (sole)
energy regularity of initial data, we can
prove (see Theorem \ref{teo:energy})
that $\teta(t)$ becomes instantaneously
bounded from below, namely, for $t>0$ we have
$\teta(t)\ge c(t)>0$, where the constant $c(t)$
depends only on the initial energy and may vanish
as $t\to 0$. If, in addition, $\teta_0$ is
(in space dimension $d=3$) in some $L^p$-space
with $p>3(=p_c)$, then (see Theorem \ref{teo:smoothing})
we prove an analogous bound from above
($\teta(t)\le C(t)<+\infty$ for $t>0$).
These bounds also entail further regularity 
properties of the solutions, which hold uniformly 
for large values of the time variable. Indeed,
$c(t)$ does not vanish, and $C(t)$ does not explode,
for $t\nearrow \infty$.
The main tool in our regularity
proof is a suitably modified Moser iteration scheme.
Note that the required extra regularity of the initial
datum is in complete agreement with the above
discussion on the case $\Omega = \mathbb{R}^d$
and with the explicit solution \eqref{explicit}.
However, since the choice of dynamic boundary
conditions \eqref{dyn:in} precludes the occurrence
both of immediate and of finite time extinction, it remains
an open challenging question to understand
whether the restriction $p > p_c$ is still optimal
in the present setting.
%
%In three dimensions the reason for considering only
%$\teta_0\in L^p$ with $p>3$ is not only technical. In fact,
%for $\Omega=\mathbb{R}^3$ one can consider the following (see \cite{Va06})
%%
%\begin{equation} \label{explicit}
%  \Theta(x,t)=2\big((T-t)^{+}\vert x\vert^{-2}\big)^{1/2},
%\end{equation}
%%
%which is a solution of the ultra fast diffusion equation
%\eqref{calorein} with $f=0$ and
%$\Theta_0(x)=2 T^{1/2}\vert x\vert^{-1}\in L^{p}_{\loc}(\mathbb{R}^3)$
%for any $p<3$. In particular, $\Theta$ is not a bounded function and
%$\Theta(\cdot,t)\notin L^{p}_{loc}(\mathbb{R}^3)$ for $p\ge 3$, at
%least until it vanishes for $t\ge T$, with $T>0$ arbitrary. As
%a consequence, we can not expect any improvement in the summability
%for initial conditions in $L^p$ with $p<3$. What happens for $p=3$??
%It is worth noting that, if we just have $\teta_0\in L^3$
%then we can exhibit (for the problem settled
%in the whole space $\RR^3$) a self-similar solution
%that does not satisfy the upper bound. This suggests
%that $3$ may play the role of a critical exponent
%also for the problem settled on a bounded domain~$\Omega$.
%
Mathematically, the main difficulties we encountered
in the analysis of this problem come from the choice of the
boundary conditions \eqref{dyn:in} which
does not allow us to perform some otherwise standard
a priori estimates. Actually, we have a sort of
``asymmetry'' of diffusion effects
between the equation in the interior domain (where
the Laplacian acts on $1/\teta$) and that on the boundary
(where the Laplace-Beltrami operator acts on (the trace of) $\teta$).
Moreover, as the precise statements of our results show,
in some situations we cannot allow $\alpha$, or $\beta$, to be $0$.
Note, however, that the situation when both $\alpha$ and $\beta$ vanish
is always permitted because it corresponds to the simpler case
of no-flux conditions, already studied in the literature. 

A natural application of equation \eqref{calorein}
comes from the so-called phase change models
of Penrose-Fife type \cite{PF}. In this physical context,
the unknown $\teta$ represents the absolute temperature
of a material liable to a phase transition,
while the source $f$ can also take into account
the effects of the phase variable on the temperature
evolution. More precisely, in the Penrose-Fife model,
equation \eqref{calorein} is coupled with
a parabolic equation of Allen-Cahn or Cahn-Hilliard type describing
the evolution of the phase variable $\chi$
(see, e.g., \cite{CGRS}, \cite{CL}, \cite{CLS},
\cite{Sch}, \cite{SSZ1}). In this framework,
the choice of considering a zero-mean-valued forcing function
$f$ in \eqref{calorein} can be motivated by the need
of replicating in our situation some inner cancellation
effects that appear in the energy estimates for the full model.
%
%The usual choice of boundary conditions for the Penrose-Fife model is that
%of Neumann boundary conditions for the phase variable and (non homogeneous)
%Dirichlet or Neumann or even of third type
%for the temperature $\teta$. With this respect, the choice of the dynamic
%boundary conditions \eqref{dyn:in} is a novelty.
%
An application of the present results to the
Penrose-Fife system with boundary conditions
of the type \eqref{dyn:in} (which has never
been studied in the literature, at least up to our
knowledge) will be given in a forthcoming paper.
Besides the Penrose-Fife model,
equation \eqref{calorein} (or, more generally, equation
\eqref{eq:diffusion} with $m < 0$) comes naturally into play
in other physical contexts (see \cite{BK} and \cite{COR}).
For example, it appears \cite{deGe} in the study of the
long-range Van der Waals interactions in thin films that
diffuse on a solid surface.

The plan of the paper is as follows.
In the next Section~\ref{sec:main},
we present our assumptions and
state a rigorous definition of weak solution.
Section~\ref{sec:moser} is devoted to proving the
main technical lemmas constituting the core of our
existence proof. In particular, using Moser iteration arguments,
we show that, under suitable conditions on data,
both $\teta$ and $\teta^{-1}$ satisfy
instantaneous regularization properties.
Then, in Section~\ref{sec:smooth} we prove that
(global) existence holds for smooth and bounded initial data.
Finally, in Section~\ref{sec:final} we prove
our main results. Namely, we show that existence
holds for initial data enjoying the sole ``energy''
regularity. Moreover, we investigate regularization properties,
and (under additional conditions) uniqueness of weak solutions.

%%%%%%%%%%%%%%%%%%%%%%%%%%%%%%%%%%%%%%%%%%%%%%%%%%%%%%%%%%%%%%%%%%%
%%%%%%%%%%%%%%%%%%%%%%%%%%%%%%%%%%%%%%%%%%%%%%%%%%%%%%%%%%%%%%%%%%%

\section{Notation and hypotheses}
\label{sec:main}

Let $\Omega$ be a smooth bounded domain of $\RR^3$
(of course, everything could be easily extended
to the one and two dimensional cases, where, actually,
better results are expected to hold).
Let also $|\Omega|=1$ so that
$\| v \|_{L^p(\Omega)}\le \| v \|_{L^q(\Omega)}$
for all $1\le p \le q \le +\infty$,
$v\in L^q(\Omega)$.
Let $H:=L^2(\Omega)$,
endowed with the standard scalar product $(\cdot,\cdot)$
and norm $\| \cdot \|$. Let also $V:=H^1(\Omega)$.
We note by $\|\cdot\|_X$
the norm in the generic Banach space $X$.
and by $\langle\cdot,\cdot\rangle_X$
the duality between $X'$ and $X$.
We will also write
$\| \cdot \|_p$ for $\| \cdot \|_{L^p(\Omega)}$
and $\| \cdot \|_{p,\Gamma}$ for $\| \cdot \|_{L^p(\Gamma)}$,
for simplicity.
Next, we set $H\Ga:=L^2(\Gamma)$
and $V\Ga:=H^1(\Gamma)$ and denote by
$(\cdot,\cdot)\Ga$ the scalar product in $H\Ga$,
by $\|\cdot\|\Ga$ the corresponding norm,
and by $\langle\cdot,\cdot\rangle\Ga$
the duality between $V'\Ga$ and $V\Ga$.
We also denote by $\nabla\Ga$ the tangential gradient
on $\Gamma$ and by $\Delta\Ga$ the Laplace-Beltrami
operator. We can thus set
\begin{equation}\label{defiHVGa}
 \calH:=H\times H_\Gamma \qquand
   \calV:=\big\{z\in V:~z|\Ga \in V\Ga\big\}.
\end{equation}
The spaces $\calH$ and $\calV$ are endowed
with the natural norms, respectively 
denoted by $\| \cdot \|_{\calH}$ and 
$\| \cdot \|_{\calV}$. For instance,
for $h\in \calH$, we may set 
$\| h \|^2_{\calH}:=\| h \|^2 + \| h \|^2\Ga$, 
whereas for $v\in \calV$ we put
$\| v \|^2_{\calV}:=\| v \|^2_{V} + \| \nabla\Ga v \|^2\Ga$. 
Unless specified otherwise, in the sequel
we shall make the following convention:
when we write $h\in\calH$, $h$ will be interpreted
as a pair of functions belonging, respectively,
to $H$ and to $H\Ga$, and
both denoted by the same letter.
On the other hand, when we consider $v\in \calV$
(or even $v\in V$), the symbol $v$ will
be intended, depending on the context,
either as a function defined on $\Omega$,
or as a pair formed by a function of $\Omega$
and its trace on $\Gamma$.

For any function, or functional $z$, defined on
$\Omega$, we can then set
\begin{equation}\label{defim}
  m\OO(z):=\frac1{|\Omega|}\io z
   = \io z,
\end{equation}
where the integral is substituted with
the duality $\langle z,1 \rangle$ in case, e.g.,
$z\in V'$. Given $\alpha\ge 0$, we also define
the measure $\dm$, given by
\begin{equation}\label{defidm}
  \ibaro v \dm:= \io v
   + \alpha \iga v,
\end{equation}
where $v$ represents a generic
function in $L^1(\Omega)\times L^1(\Gamma)$.
Here and below, integrals over $\Gamma$ are to
be intended with respect to the standard surface
measure. With some abuse of notation, we will
also write
\begin{equation}\label{defimm}
  m(v):= \frac1{|\Omega| + \alpha |\Gamma|}
   \ibaro v \dm = \frac1{1 + \alpha |\Gamma|}
   \ibaro v \dm,
\end{equation}
i.e., the ``mean value'' of $v$ w.r.t.~the measure
$\dm$. Here $|\Gamma|$ represents the surface measure
of $\Gamma$. In case $\alpha=0$, it is intended that
$m(v)=m\OO(v)$.
For $p\in[1,\infty)$ and $X$ a Banach space,
we introduce the space
\begin{equation}\label{defT1}
  \calT^p(0,+\infty;X)
   := \Big\{g\,\in\,L^p_{\loc}(0,+\infty;X)\,:~
   \sup_{t\in[0,+\infty)} \| g \|_{L^p(t,t+1;X)} < +\infty \Big\},
\end{equation}
which is a Banach space as it is endowed with the graph norm.
Assuming $f:(0,+\infty)\times \Omega\to \RR$ be
a suitable source term and letting $\alpha,\beta\ge 0$,
we can introduce the system
\begin{align}\label{calore}
  & \teta_t - \Delta u = f,
   \quad u = -\frac1\teta,
   \quext{in }\,(0,+\infty)\times\Omega,\\
 \label{calorega}
  & \alpha \eta_t - \beta \Delta\Ga \eta = - \dn u,
   \quext{on }\,(0,+\infty)\times\Gamma,\\
 \label{iniz}
  & \teta|_{t=0}=\teta_0
   \quext{in }\,\Omega,\\
 \label{inizga}
  & \alpha \eta|_{t=0}=\alpha \eta_0
   \quext{on }\,\Gamma.
\end{align}
Our basic assumptions on the initial data
are the following:
\begin{align}\label{hpteta0w}
  & \teta_0 \in L^1(\Omega),\quad
    \log\teta_0 \in L^1(\Omega),\\
 \label{hpeta0w}
  & \alpha\eta_0 \in L^1(\Gamma),\quad
     \alpha\log\eta_0 \in L^1(\Gamma).
\end{align}
This very natural condition corresponds to
asking that the initial data have {\sl finite energy},
where the energy functional $\calE$ is defined as
\begin{equation}\label{defiE}
  \calE(\teta,\eta):=
   \io\big(\teta-\log\teta\big)
   +\alpha\iga \big(\eta-\log\eta\big).
\end{equation}
In the sequel, we will simply write, with
some abuse of language, $\calE(t)$,
in place of $\calE(\teta(t),\eta(t))$.
Moreover, we will note as $\EE_0$ the energy of 
initial data, namely
\begin{equation}\label{defiE0}
  \EE_0:=\calE(\teta_0,\eta_0).
\end{equation}
The source term is assumed to satisfy
for some given $\epsilon\in(0,1)$
\begin{equation}\label{hypfw}
  f\in L^2(0,+\infty;L^{6/5}(\Omega))
   \cap \mathcal{T}^2(0,+\infty;L^{3+\epsilon}(\Omega)), \qquad
   m\OO(f)(t) = 0,\quext{for a.a.~}\,t\in(0,+\infty).
\end{equation}
The first regularity condition seems to be
necessary for controlling the energy {\sl uniformly}\/ in time. 
The second condition, where we require more summability in space,
but allow for a {\sl locally uniform}\/ (rather than uniform) 
summability in time, is used for the purpose of Moser iterations. 
The assumption $m\OO(f)(t) = 0$ is taken 
because we need a uniform estimate for the spatial mean 
both of $\teta$ and of $u$, which may not hold
for general $f$. Actually, integrating \eqref{calore}
in $\Omega$, \eqref{calorega} on $\Gamma$ and taking the sum,
one can see that the ``total mass'' $m(\teta)$ is conserved in time.
It is also worth observing that it may be possible to
extend our results by considering more general conditions on
data (and, particularly, on~$f$). However we believe
that the assumptions provided above are very natural,
especially in relation with the energy estimate, and 
may allow us to give simpler proofs.
 
In the sequel we will note as {\sl Problem~(P)}
the the initial-value problem
for system \eqref{calorein}-\eqref{dyn:in}. In particular,
we can introduce a suitable concept of weak solution
as follows: 
\bede\label{def:weak}
 A (global) weak solution (or ``energy solution'')
 to\/ {\rm Problem~(P)} is a triplet $(\teta,\eta,u)$ 
 satisfying, for all $T>0$, the regularity properties
 \begin{align}\label{regotetaw}
   & \teta\in C^0([0,T];L^1(\Omega)), \qquad
    \alpha \eta\in C^0([0,T];L^1(\Gamma)),\\
  \label{signteta}
   & \teta>0 \quext{a.e.~in }\,(0,T)\times\Omega, \qquad
    \eta>0 \quext{a.e.~on }\,(0,T)\times\Gamma,  \\
  \label{regouw}
   & u\in L^2(0,T;V), \qquad u=-1/\teta, \quext{a.e.~in }\,(0,T)\times\Omega,\\
  \label{regoetaw}
   & \beta \log \eta\in L^2(0,T;V\Ga),\\
  \label{trace}
   & \eta=-1/u\Ga, \quext{a.e.~on }\,(0,T)\times \Gamma,
 \end{align}
 and fulfilling, for any test function
 \begin{equation}\label{regotest}
   \xi \in C^1([0,T];C^0(\barO))\cap
    C^0([0,T];C^2(\barO))
 \end{equation}
 and for all times $t\in[0,T]$, the relation
 \begin{align}\no
   & \io \teta(t) \xi(t) + \alpha \iga \eta(t) \xi(t)
    + \ito \nabla u \cdot \nabla \xi
    - \beta \itga \eta \Delta\Ga \xi \\
   \label{caloreveryw}
   & \mbox{}~~~~~~~~~~
    = \ito \teta \xi_t
    + \alpha \itga \eta \xi_t
    + \ito f \xi
    + \io \teta_0 \xi(0)
    + \alpha \iga \eta_0 \xi(0).
 \end{align}
\edde
\beos\label{tetaeta}
 It is worth giving some explanation of relation \eqref{trace}.
 There, $u\Ga$ denotes the trace of $u$ on $\Gamma$, which exists
 for almost every value of the time variable thanks to
 the first~\eqref{regouw}. More precisely, we have
 $u\Ga\in L^2(0,T;H^{1/2}(\Gamma))$. Actually, for weak
 solutions we cannot simply write $\eta=\teta\Ga$ since
 the trace of $\teta$ does not necessarily exist. On the other hand,
 when considering smoother solutions (e.g., in the approximation
 detailed in Section~\ref{sec:smooth} below)
 it will happen that $\teta$ is more regular
 so that $\eta$ can be intended, in fact, as the trace 
 of~$\teta$.
\eddos
\beos\label{sualfabeta}
 In the case when $\alpha=0$ and $\beta>0$, the last integral on the
 first row of \eqref{caloreveryw} may make no sense since 
 \eqref{regotetaw}-\eqref{trace} do not guarantee
 any summability property for $\eta$, whereas 
 when $\alpha>0$ we can get help from the latter \eqref{regotetaw}.
 For this reason we will be able to consider 
 the case $\alpha=0$, $\beta>0$ only for more regular solutions
 (and, correspondingly, under more restrictive
 conditions on data).
\eddos
%

%%%%%%%%%%%%%%%%%%%%%%%%%%%%%%%%%%%%%%%%%%%%%%%%%%%%%%%%%%%%%%%%%%%%%%%%%%%%%%%
%%%%%%%%%%%%%%%%%%%%%%%%%%%%%%%%%%%%%%%%%%%%%%%%%%%%%%%%%%%%%%%%%%%%%%%%%%%%%%%

\section{Main technical lemmas}
\label{sec:moser}

In this section we prove some regularization estimates
holding for sufficiently smooth solutions of Problem~(P).
In this procedure, a sufficient regularity will always
be assumed in such a way that all computations we perform
make sense. For this reason, we will indicate simply by
$\teta$, rather than by $\eta$, the boundary temperature
(indeed, the ``bulk'' $\teta$ has a trace at this regularity
level). As a rule, we will use the letters $c$ and $\kappa$
to denote generic positive constants, depending only
on the set $\Omega$ and on the parameters $\alpha$ and $\beta$,
with $\kappa$ used in estimates from below.
The dependence, or independence, of $\kappa$ and 
$c$, with respect to time will be pointed out on occurrence. 
The values of $c$, $\kappa$ are allowed to vary
from time to time.
Finally, $Q$ will denote a generic computable
function, increasingly monotone with respect to each of its
parameters, taking values in $[0,+\infty)$.
%

%%%%%%%%%%%%%%%%%%%%%%%%%%%%%%%%%%%%%%%%%%%%%%%%%%%%%%%%%%%%%%%%%%%

\subsection{Regularization of $u$}
\label{subsec:moseru}

\begin{lemm}[Moser iterations for $u$]%
 \begin{sl}\label{moser-u}
 Let $S\in\RR$, $T\ge S+2$,  
 $\alpha\ge 0$, $\beta\ge 0$, and
 let, for some $\epsilon\in(0,1)$,
 \begin{equation}\label{inpiuf}
   f \in L^2(S,T;L^{3+\epsilon}(\Omega)), \qquad
    \| f \|_{L^2(S,T;L^{3+\epsilon}(\Omega))} =: F.
 \end{equation}
 Let $(u,\teta)$ be a couple of sufficiently
 smooth functions solving, in a suitable sense,
 the system
 \begin{align}\label{eq-abstr}
   & \teta_t - \Delta u = f, \quad \teta=-1/u,
    \quext{in }\,\Omega,\\
  \label{bc-abstr}
   & \alpha \teta_t - \beta \Delta\Ga \teta = -\dn u,
   \quext{on }\,\Gamma,
 \end{align}
 over the time interval $(S,T)$.
 Moreover, let us assume that
 \begin{align}\label{inpiuu}
  & u \in L^2(S,T;V), \qquad
    \| u \|_{L^2(S,T;V)} =: M,\\
  \label{uSL1}
    & u(S)\in L^1(\Omega),~~\alpha u(S)\in L^1(\Gamma), \qquad
      \| u(S) \|_{L^1(\Omega)}
     + \alpha \| u(S) \|_{L^1(\Gamma)} =: U.
 \end{align}
 Then, for any $\tau\in (0,1)$, we have
 \begin{equation}\label{st-abstr-u}
   \| u \|_{L^\infty((S+\tau,T)\times \Omega)} \le
    Q\big( F, M, U, \tau^{-1}, T-S \big).
 \end{equation}
\enle
\begin{proof}
We will just consider the case when $\beta=0$, which is more difficult
since we cannot get any help from the boundary diffusion term.
In this proof, the constant(s) $c$ 
are allowed to depend (monotonically) on the length of the time interval $(S,T)$, 
which is finite and assigned. Indeed, in the sequel this lemma will be applied on time
intervals of fixed length. 
Setting $z:=-u\ge 0$, we will also assume, for the sake of simplicity,
that $z\ge1$ almost everywhere. Indeed, if that does not hold, then it
is easy to check that in the estimates below we can simply replace
$z$ with $z=\max\left\{z,1\right\}$. That said, we
test \eqref{eq-abstr} by $z^{p+1}$, where $p\ge 1$
will be specified later. This gives
\begin{equation}\label{lemma1-1}
  \frac1p \ddt \big( \| z \|_p^p + \alpha \| z \|_{p,\Gamma}^p \big)
   + \frac{4(p+1)}{(p+2)^2} \big\| \nabla z^{\frac{p+2}2} \big\|^2
  %
  % + \frac{4\beta(p+1)}{p^2} \big\| \nabla\Ga z^{\frac{p}2} \big\|\Ga^2
  %
   \le \io |f| z^{p+1}.
\end{equation}
In order to recover the full $V$-norm on the \lhs,
we can multiply \eqref{lemma1-1} by $p$ and
then add to both hand sides the quantity
\begin{equation}\label{lemma1-3}
  \big\| z^{\frac{p+2}2} \big\|_1^2
   = \| z \|_{\frac{p+2}2}^{p+2}.
\end{equation}
Hence, using continuity of the embedding
$V\subset L^6(\Omega)$ and of the trace operator
from $V$ to $L^4(\Gamma)$, we get on the \lhs\ of \eqref{lemma1-1}
a quantity $\calI$ such that
\begin{equation}\label{lemma1-a1}
  \calI := \kappa \big\| \nabla z^{\frac{p+2}2} \big\|^2
   + \| z \|_{\frac{p+2}2}^{p+2}
   \ge \kappa \| z \|_{3(p+2)}^{p+2}
   + \kappa \| z \|_{2(p+2),\Gamma}^{p+2}.
\end{equation}
We now choose $p=1$ for the first iteration. Then,
we can estimate the \rhs\ of \eqref{lemma1-3}
as
\begin{equation}\label{lemma1-3b}
  \| z \|_{\frac{p+2}2}^{p+2}
   \stackrel{p=1}{=} \| z \|_{\frac{3}2}^{3}
   \le c \| z \|_6^2 \| z \|_1.
\end{equation}
On the other hand, still for $p=1$, we can write
\begin{equation}\label{lemma1-2-1}
  \io |f| z^{2}
   \le \| f \|_{3} \| z^{1/2} \| \| z^{3/2} \|_6
   \le \sigma \| z \|_9^3
   + c_\sigma \| f \|_3^2 \| z \|_1,
\end{equation}
whence, taking $\sigma$ small enough
and collecting \eqref{lemma1-3}-\eqref{lemma1-2-1},
\eqref{lemma1-1} becomes
\begin{equation}\label{lemma1-1b}
  \ddt \big( \| z \|_1 + \alpha \| z \|_{1,\Gamma} \big)
   + \kappa \| z \|_9^3
   + \kappa \| z \|_{6,\Gamma}^3
   \le c \big ( \| z \|_6^2 + \| f \|_3^2\big) \| z \|_1.
\end{equation}
Thus, integrating over $(S,T)$ and using Gronwall's
lemma, we readily arrive at
\begin{equation}\label{iter-1}
  \| z \|_{L^\infty(S,T;L^1(\Omega))}
   + \alpha \| z \|_{L^\infty(S,T;L^1(\Gamma))}
   + \| z \|_{L^3(S,T;L^9(\Omega))}
   + \| z \|_{L^3(S,T;L^6(\Gamma))}
 %  + \beta \| u \|_{L^3((S,S+2;L^6(\Gamma))}
  \le Q\big( F, M, U \big),
\end{equation}
where the function $Q$ has no
explicit dependence on $T$ at this level. This relation
is the starting point for the subsequent iterations.

To proceed, we set $r:=\frac{3+\epsilon}{2+\epsilon}$
to be the conjugate exponent of $3+\epsilon$.
Rewriting \eqref{lemma1-1} (multiplied by $p$)
for a suitable new choice
of $p$, we can now estimate the \rhs\ as follows:
\begin{equation}\label{lemma1-2}
  p\io |f| z^{p+1}
   \le p \| f \|_{3+\epsilon} \big\| z^{p+1} \big\|_r
   = p \| f \|_{3+\epsilon} \| z \|_{r(p+1)}^{p+1}.
\end{equation}
Then, adding again the term $\| z \|_{(p+2)/2}^{p+2}$
to both hand sides and integrating the result from
some $\ttt\ge S$ (to be chosen later) to $T$,
we arrive at
\begin{align} \no
  & \| z \|_{L^{\infty}(\ttt,T;L^p(\Omega))}^p
   + \alpha \| z \|_{L^{\infty}(\ttt,T;L^p(\Gamma))}^p
   + \| z \|_{L^{p+2}(\ttt,T;L^{3(p+2)}(\Omega))}^{p+2}
   + \| z \|_{L^{p+2}(\ttt,T;L^{2(p+2)}(\Gamma))}^{p+2}\\
 %   + \beta\, \text{(pezzo di bordo)}\\
 \no
  & \mbox{}~~~
   \le c \| z(\ttt) \|_{p}^p
   + c \alpha \| z(\ttt) \|_{p,\Gamma}^p
   + c p F \| z \|_{L^{2(p+1)}(\ttt,T;L^{r(p+1)}(\Omega))}^{p+1}
   + c \| z \|_{L^{p+2}(\ttt,T;L^{\frac{p+2}2}(\Omega))}^{p+2}\\
 \label{lemma1-4}
  & \mbox{}~~~
   \le c \| z(\ttt) \|_{p}^p
   + c \alpha \| z(\ttt) \|_{p,\Gamma}^p
   + c ( p F + 1)  \| z \|_{L^{2(p+1)}(\ttt,T;L^{r(p+1)}(\Omega))}^{p+2},
\end{align}
where we used, in particular, that $z\ge 1$. 
%
% We may note that the last constant $c$ is actually dependent (monotonically) 
%also on the length of the time interval $(S,T)$, which is however assigned 
%a priori. Actually, in the sequel this lemma will be applied on time
%intervals of fixed length. 
%
Next, we define, for $i\ge 0$ and $\tau_i$ to be chosen later,
\begin{equation}\label{Jp}
  J_i^{p_i} := \| z \|^{p_i}_{L^\infty(\tau_i,T;L^{p_i}(\Omega))}
   + \| z \|^{p_i+2 }_{L^{p_i+2}(\tau_i,T,L^{3p_i+6}(\Omega))}
   + \| z \|_{L^{p_i+2}(\ttt,T;L^{2(p_i+2)}(\Gamma))}^{p_i + 2 }.
\end{equation}
Then, using that $z\ge 1$ and that the integration domain
has measure greater than $1$, 
by elementary interpolation we obtain
\begin{align}\nonumber
  \| z \|_{L^{2(p_{i+1}+1)}(\tau_i,T;L^{r(p_{i+1}+1)}(\Omega))}
   & \le \| z \|_{L^\infty(\tau_i,T;L^{p_i}(\Omega))}^{\rho}
    \| z \|_{L^{p_i+2}(\tau_i,T,L^{3p_i+6}(\Omega))}^{1-\rho}\\
 \nonumber
  & \le \rho \| z \|_{L^\infty(\tau_i,T;L^{p_i}(\Omega))}
    + (1-\rho) \| z \|_{L^{p_i+2}(\tau_i,T,L^{3p_i+6}(\Omega))}\\
 \label{Jpnew}
  & \le \rho \| z \|_{L^\infty(\tau_i,T;L^{p_i}(\Omega))}
    + (1-\rho) \| z \|^{\frac{p_i+2}{p_i}}_{L^{p_i+2}(\tau_i,T,L^{3p_i+6}(\Omega))}, 
\end{align}
where the index $p_{i+1}$ and the interpolation exponent $\rho=\rho(i)$
are given by the system
\begin{equation} \label{interpol12}
  \displaystyle
   \begin{cases}
     \frac{1-\rho}{p_i+2} = \frac{1}{2(p_{i+1}+1)},\\
      \frac{\rho}{p_i} + \frac{1-\rho}{3(p_i+2)}\,=\,\frac{1}{r(p_{i+1}+1)}.
   \end{cases}
\end{equation}
Dividing the second equation in \eqref{interpol12}
by the first one, we have
\begin{equation} \label{interpol13}
  \Big(\frac{\rho}{p_i} + \frac{1-\rho}{3(p_i+2)} \Big)\frac{p_i+2}{1-\rho}=\frac{2}{r},
\end{equation}
whence
\begin{equation} \label{interpol14}
  \frac{\rho}{p_i}\frac{p_i+2}{1-\rho}=\frac{2}{r}-\frac{1}{3}=:K_\epsilon,
\end{equation}
and it is easy to compute
\begin{equation} \label{interpol15}
  K_\epsilon=\frac{9+5\epsilon}{9+3\epsilon}.
\end{equation}
From \eqref{interpol14} we also infer
\begin{equation} \label{interpol16}
  \rho = \frac{ \frac{p_i}{p_i+2} K\ee }{ 1 + \frac{p_i}{p_i+2} K\ee }
   \in (0,1), \quext{provided that }\, p_i\ge 1.
\end{equation}
Being
\begin{equation} \label{interpol17}
  1 - \rho = \frac{ 1 }{ 1 + \frac{p_i}{p_i+2} K\ee },
\end{equation}
we then obtain from the first \eqref{interpol12}
\begin{align} \no
  p_{i+1} & = \frac12 \frac{p_i+2}{1-\rho} - 1
    = \frac12 \Big(1 + \frac{p_i}{p_i+2} K\ee \Big) (p_i+2) - 1\\
  \label{interpol18}
  & = \frac{ K\ee + 1 }2 p_i = \frac{9+4\epsilon}{9+3\epsilon} p_i =: H p_i,
\end{align}
where, obviously, $H=H(\epsilon)>1$ whenever $\epsilon>0$.

Taking the $p_i$-th power of \eqref{Jpnew},
by convexity we obtain 
\begin{equation}\label{Jpubis}
  J_i \ge \| z \|_{L^{2(p_{i+1}+1)}(\tau_i,T;L^{r(p_{i+1}+1)}(\Omega))}.
\end{equation}
We can now start the iteration argument. Let $p_0=1$ and
inductively define, for $i \ge 1$, $p_{i+1} := H p_i = H^{i+1}$.
Moreover, let (for instance), for $i\ge 1$,
\begin{equation} \label{sigmai}
  \ttt_i := \tau \frac3{\pi^2i^2}, \quext{so that }\,
   \sum_{i=1}^\infty \ttt_i = \frac{\tau}2.
\end{equation}
Let us now rewrite \eqref{lemma1-4}
by taking $p=p_{i+1}$.
Setting also, for brevity, $J_i:=J_{p_i}$
and writing $\tau_{i+1}$ in place of $\ttt$
(to be chosen below), we
then obtain, thanks also to \eqref{Jpubis},
\begin{equation}
  J_{i+1}^{p_{i+1}}
   + \alpha \| z \|^{p_{i+1}}_{L^\infty(\tau_{i+1},T;L^{p_{i+1}}(\Gamma))}
%   + \| z \|^{p_{i+1}+2}_{L^{p_{i+1}+2}(\tau_{i+1},T;L^{2(p_{i+1}+2)}(\Gamma))} \\
 \label{moser10}
 %  & \mbox{}~~~~~
    \le c \| z(\tau_{i+1}) \|_{p_{i+1}}^{p_{i+1}}
    + c \alpha \| z(\tau_{i+1}) \|_{p_{i+1},\Gamma}^{p_{i+1}}
    + c \big( p_{i+1} F + 1\big) J_i^{p_{i+1}+2}.
\end{equation}
Let us now make precise the choice of $\tau_i$.
First of all, we set $\tau_0:= S$. Then, 
by induction, given $\tau_i$, we take
$\tau_{i+1} \in (\tau_i,\tau_i+\ttt_{i+1})$
such that
\begin{align} \no
  & \| z(\tau_{i+1}) \|_{p_{i+1}}^{p_{i}}
  + \alpha \| z(\tau_{i+1}) \|_{p_{i+1},\Gamma}^{p_{i}}
  \le \frac1{\ttt_{i+1}}\int_{\tau_i}^{\tau_i+\ttt_{i+1}}
   \Big( \| z(t) \|_{p_{i+1}}^{p_{i}}
    + \alpha \| z(t) \|_{p_{i+1},\Gamma}^{p_{i}} \Big)\,\dit\\
  \label{moser10b}
  & \mbox{}~~~~~~~~~~
  \le c \frac{i^2}{\tau} \int_{\tau_i}^{\tau_i+\ttt_{i+1}}
    \Big( \| z(t) \|_{3(p_i+2)}^{p_{i}}
    + \alpha \| z(t) \|_{2(p_i+2),\Gamma}^{p_{i}} \Big)\,\dit
  \le c \frac{i^2}{\tau} J_i^{p_i},
\end{align}
where we used that $p_{i+1}\le 2(p_i+2)$.
Thus, \eqref{moser10}-\eqref{moser10b} give
\begin{equation}
  J_{i+1}^{p_{i+1}}
   + \alpha \| z \|^{p_{i+1}}_{L^\infty(\tau_{i+1},T;L^{p_{i+1}}(\Gamma))}
%   + \| z \|^{p_{i+1}+2}_{L^{p_{i+1}+2}(\tau_{i+1},T;L^{2(p_{i+1}+2)}(\Gamma))} \\
 \label{moser11}
%  \mbox{}~~~~~~~~~
    \le c i^{2H}\tau^{-H} J_i^{p_{i+1}}
    + c \big( p_{i+1} F + 1\big) J_i^{p_{i+1}+2}.
\end{equation}
In particular, we obtain
\begin{equation} \label{moser12}
  J_{i+1}^{H^{i+1}} \le c \Big( i^{2H}\tau^{-H} + H^{i+1} F + 1 \Big)
   J_i^{H^{i+1}+2}.
\end{equation}
Thus, setting $\eta_i:=(H^i+2)/H^i$, we have
\begin{equation} \label{moser13}
  J_{i+1} \le B_i^{H^{-(i+1)}}
   J_i^{\eta_{i+1}}, \quext{where }\,
   B_i:= c \big( i^{2H}\tau^{-H} + H^{i+1} F + 1 \big).
\end{equation}
Hence, letting $L_i:=\log J_i$ and
$\zeta_i:=\log B_i$, we can rewrite
\eqref{moser13}, for $i$ large enough, as
\begin{align} \no
  L_{i+1} & \le H^{-(i+1)} \zeta_i + \eta_{i+1} L_i
   \le H^{-(i+1)} \zeta_i + \eta_{i+1} \big( H^{-i} \zeta_{i-1} + \eta_i L_{i-1} \big) \\
 \no
  & \le H^{-(i+1)} \zeta_i + \eta_{i+1} H^{-i} \zeta_{i-1} + \eta_{i+1} \eta_i
   \big( H^{-(i-1)} \zeta_{i-2} + \eta_{i-1} L_{i-2} \big) \\
% \no
%  & \le H^{-(i+1)} \zeta_i + \eta_{i+1} H^{-i} \zeta_{i-1}
%  + \eta_{i+1} \eta_i H^{-(i-1)} \zeta_{i-2}
%  + \eta_{i+1} \eta_i \eta_{i-1} \big( H^{-(i-2)} \zeta_{i-3} + \eta_{i-2} L_{i-3} \big) \\
 \label{moser14}
  & \le \dots
   \le L_0 \prod_{k=1}^{i+1} \eta_k
    + \sum_{k=1}^{i+1} \Big( \zeta_{k-1} H^{-k} \prod_{j=k+1}^{i+1} \eta_j \Big),
\end{align}
where the last productory is understood to be $1$ for $k=i+1$. 
Then, letting $i\nearrow \infty$, it is not difficult to 
obtain \eqref{st-abstr-u}, which concludes the proof
(see also \cite{Sch} for more details).
\end{proof}
\noindent%
If $u$ is bounded at the initial time, we can avoid all
complications connected with the choice of the sequence $\ttt_i$.
Indeed, a straightforward modification of the above proof
permits to show the following
\beco\label{moser-u-2}
 Let $S\in\RR$, $T\ge S+2$, $\alpha\ge 0$, $\beta\ge 0$, and
 let, for some $\epsilon\in(0,1)$, \eqref{inpiuf} hold.
 Let, as before, $(u,\teta)$ be a couple of sufficiently
 smooth functions solving, in a suitable sense,
 system\/ \eqref{eq-abstr}-\eqref{bc-abstr}
 over the time interval $(S,T)$.
 Moreover, let us assume \eqref{inpiuu}
 together with
 \begin{equation} \label{uSLinfty}
    u(S)\in L^\infty(\Omega),~~\alpha u(S)\in L^\infty(\Gamma), \qquad
      \| u(S) \|_{L^\infty(\Omega)}
     + \alpha \| u(S) \|_{L^\infty(\Gamma)} =: U.
 \end{equation}
 Then, we have 
 \begin{equation}\label{st-abstr-u-2}
   \| u \|_{L^\infty((S,T)\times \Omega)} \le
   Q\big( F, M, U, T-S \big). 
 \end{equation}
\enco
%
%

%%%%%%%%%%%%%%%%%%%%%%%%%%%%%%%%%%%%%%%%%%%%%%%%%%%%%%%%%%%%%%%%%%%

\subsection{Regularization of $\teta$}
\label{subsec:moserteta}

\begin{lemm}[uniform stability of $L^p$ norm]%
 \begin{sl}\label{moser-teta-1}
 Let $S\in \RR$, $T\ge S+2$, $\alpha\ge 0$, $\beta\ge 0$.
 Let $(u,\teta)$ be a couple of sufficiently smooth
 functions solving in a suitable sense system\/
 \eqref{eq-abstr}-\eqref{bc-abstr} over $(S,T)$.
 Let also the following properties hold:
 \begin{align}\label{tetaf-1}
   & f \in L^2(S,T;L^3(\Omega)), \qquad
    \| f \|_{L^2(S,T;L^3(\Omega))} =: F,\\
   \label{tetaf-2}
   & \teta \in L^\infty(S,T;L^1(\Omega)), \qquad
   \| \teta \|_{L^\infty(S,T;L^1(\Omega))} =: L,\\
  \label{tetaz}
   & \teta(S)\in L^p(\Omega),~~\alpha \teta(S)\in L^p(\Gamma), \qquad
    \| \teta(S) \|_{L^{p}(\Omega)}
     + \alpha \| \teta(S) \|_{L^{p}(\Gamma)} =: \Theta,
 \end{align}
 where the exponent $p$ is given and is assumed to satisfy $p\in [4,\infty)$
 in the case when $\alpha > 0$ and $\beta = 0$, and is assumed to satisfy
 $p\in [3,\infty)$ in all the other cases. 
 Then,
 \begin{equation}\label{st-abstr-teta0}
   \| \teta \|_{L^\infty(S,T;L^{p}(\Omega))}
   + \alpha \| \teta \|_{L^\infty(S,T;L^{p}(\Gamma))}
   \le Q( L, F, \Theta ),
 \end{equation}
 where the expression of $Q$ is\/ {\rm independent} of the final time $T$. 
\enle
\begin{proof}
Let us consider first the case when $p\le 4$.
We test \eqref{eq-abstr} by $\teta^{p-1}$. This gives
\begin{equation}\label{lemma2-11}
  \ddt \Big( \| \teta \|_p^p + \alpha \| \teta \|_{p,\Gamma}^p \Big)
   + \frac{4p(p-1)}{(p-2)^2} \big\| \nabla \teta^{\frac{p-2}2} \big\|^2
   + \beta \frac{4p(p-1)}{p^2} \big\| \nabla\Ga \teta^{\frac{p}2} \big\|^2\Ga
   \le p \io |f| \teta^{p-1}.
\end{equation}
Then, being $p\le 4$, using \eqref{tetaf-2},
we can add the inequality
\begin{equation}\label{lemma2-12}
  \big\| \teta^{\frac{p-2}2} \big\|_1^2 
   + \beta \big\| \teta \big\|_1^p 
   \le Q(L).
\end{equation}
Then, owing also to the continuity of the trace operator from 
$V$ to $L^4(\Gamma)$, we may conclude that 
\begin{align}\label{newtr1}
  & \frac{4p(p-1)}{(p-2)^2} \big\| \nabla \teta^{\frac{p-2}2} \big\|^2
   + \big\| \teta^{\frac{p-2}2} \big\|_1^2 
   \ge \kappa \| \teta \|_{3(p-2)}^{p-2} 
    + \kappa \| \teta \|_{2(p-2),\Gamma}^{p-2},\\
  \label{newtr2}
  & \beta \frac{4p(p-1)}{p^2} \big\| \nabla\Ga \teta^{\frac{p}2} \big\|^2\Ga
   + \beta \big\| \teta \big\|_1^p 
   \ge \kappa \beta \| \teta \|_{p,\Gamma}^{p}.
\end{align}
Here and below, the constants $c$ and $\kappa$ are assumed to be 
independent of $p$ and of time.
Next, we estimate the \rhs\ of \eqref{lemma2-11}
as follows:
\begin{equation}\label{lemma2-13}
  p \io |f| \teta^{p-1}
   \le c \| f \|_{3} \| \teta^{\frac{p}2} \|_2 \| \teta^{\frac{p-2}2} \|_6
    \le \sigma \| \teta \|_{3(p-2)}^{p-2}
    + c_\sigma \| f \|_{3}^2 \| \teta \|_p^p.
\end{equation}
Taking $\sigma$ small enough, we then arrive at the inequality
\begin{equation}\label{lemma2-14}
  \ddt \big( \| \teta \|_p^p
    + \alpha \| \teta \|_{p,\Gamma}^p  \big)
    + \kappa \| \teta \|_{3(p-2)}^{p-2}
    + \kappa \| \teta \|_{2(p-2),\Gamma}^{p-2}
    + \kappa \beta \| \teta \|_{p,\Gamma}^{p}  
    \le c \| f \|_{3}^2 \| \teta \|_p^p
    + Q(L).
\end{equation}
Let us now set 
\begin{equation}\label{def-y}
  y(t):= \big( \| \teta \|_p^p + \alpha \| \teta \|_{p,\Gamma}^p \big)^{1/p}. 
\end{equation}
Then, we need to distinguish between two different cases. First, let either
$\alpha = 0$ or both $\alpha$ and $\beta$ be strictly positive. Then, 
thanks to the assumption $p\ge 3$, it follows that 
$3(p-2) \ge p$. Moreover, assuming for simplicity $\teta\ge 1$, which 
is of course not restrictive, \eqref{lemma2-14} reduces to 
the differential inequality
\begin{equation}\label{lemma2-14b}
  \ddt y^p(t)
    + \kappa y^{p-2}(t)
    \le c \| f(t) \|_{3}^2 y^p(t)
    + Q(L).
\end{equation}
Then, integrating explicitly the ODE and 
applying the comparison principle,
it is not difficult to obtain \eqref{st-abstr-teta0}.
We observe in particular that
the procedure works on time intervals of any
length with no explicit dependence of $Q$ on time. In other
words, the resulting function 
$Q$ depends on $T$ only through the quantities $F$ and $L$.
Note also that the dissipative term $\kappa y^{p-2}(t)$
on the \lhs\ is essential for the purpose of getting rid
of the quantity $Q(L)$ on time intervals of arbitrary length.

Second, let us consider the case when $\alpha >0 $ and $\beta = 0$. Then,
we need to use the boundary term resulting from \eqref{newtr1}. This is
what forces us to assume $p\ge 4$ so that $2(p-2) \ge p$ (hence, in 
fact, we are taking exactly $p=4$ at least as a first iteration). 
Under these assumptions, we still deduce a differential inequality
of the form \eqref{lemma2-14b}, which is integrated exactly as before.

Finally, let us remove the restriction $p\le 4$ made at the beginning.
In that case, we need to iterate the procedure a {\sl finite}\/ number
of times. Namely, we can first take $p=4$ to
get an estimate of $\teta$ in $L^\infty(S,T;L^4(\Omega))$
(and of $\alpha\teta$ in $L^\infty(S,T;L^4(\Gamma))$).
Thanks to this new information, we can repeat the procedure by taking
higher values of $p$. We omit the details and just point out that
a similar argument will be carried out in the Lemmas~\ref{moser-teta-2},
\ref{moser-teta-2b} below.    
\end{proof}
\noindent%
\begin{lemm}[Moser iterations for $\teta$]%
 \begin{sl}\label{moser-teta-2}
 Let $S\in \RR$, $T\ge S+2$, $\alpha\ge 0$, $\beta\ge 0$.
 Moreover, if $\beta=0$ let also $\alpha=0$.
 Let $(u,\teta)$ be a couple of sufficiently smooth
 functions solving in a suitable sense system\/
 \eqref{eq-abstr}-\eqref{bc-abstr} over $(S,T)$.
 Moreover, let us assume that, for some $\epsilon\in(0,1)$,
 \begin{align}\label{tetaf2-1}
   & f \in L^2(S,T;L^{3+\epsilon}(\Omega)), \qquad
    \| f \|_{L^2(S,T;L^{3+\epsilon}(\Omega))} =: F,\\
   \label{tetaf2-2}
   & \teta \in L^\infty(S,T;L^1(\Omega)), \qquad
   \| \teta \|_{L^\infty(S,T;L^1(\Omega))} =: L,\\
  \label{tetaz2}
   & \teta(S)\in L^{3+\epsilon}(\Omega),~~
         \alpha \teta(S)\in L^{3+\epsilon}(\Gamma), \qquad
    \| \teta(S) \|_{L^{3+\epsilon}(\Omega)}
     + \alpha \| \teta(S) \|_{L^{3+\epsilon}(\Gamma)} =: \Theta.
 \end{align}
 Then, for any $\tau\in (0,1)$,
 we have
 \begin{equation}\label{st-abstr-teta}
   \| \teta \|_{L^\infty((S+\tau,T)\times \Omega)}
   + \alpha \| \teta \|_{L^\infty((S+\tau,T)\times\Gamma)}
   \le Q( F, L, \Theta, \tau^{-1}, T-S ). 
 \end{equation}
\enle
\begin{proof}
In this proof, we work on time intervals of assigned
finite length. Hence, we may assume, as in Lemma~\ref{moser-u},
explicit dependence of $c$ on the length $T-S$. 
To start with a further Moser iteration procedure,
we rewrite \eqref{lemma2-11}, multiply it by $p$,
and then add to both sides the
term $\| \teta^{\frac{p-2}2}\|_1^2$, which may be 
controlled as in \eqref{lemma2-12}. 
We then deduce
\begin{equation}\label{lemma2-21}
  \ddt \big( \| \teta \|_p^p + \alpha \| \teta \|_{p,\Gamma}^p \big)
   + \kappa \big\| \teta \big\|^{p-2}_{3(p-2)}
   + \beta \kappa \big\| \teta \big\|^{p-2}_{3(p-2),\Gamma}
   \le p \io |f| \teta^{p-1}
   + c \| \teta \|_{\frac{p-2}2}^{p-2}.
\end{equation}
Then, in the case $\alpha=\beta=0$, the proof works
similarly with \cite[Lemma~3.5]{SSZ1}, to which we
refer the reader for details. In the cases $\alpha,\beta>0$
and $\alpha=0$, $\beta>0$, the argument can be adapted
just with small variants. Note in particular that, at
the first iteration step, we choose $p=p_0=3+\epsilon$
so that the latter term in \eqref{lemma2-21}
can be controlled by using Lemma~\ref{moser-teta-1}.
\end{proof}
\noindent%
In the case when $\beta=0$ and $\alpha>0$, we can still prove
regularization of $\teta$, but the argument is a bit more delicate
since, as in the corresponding case in Lemma~\ref{moser-teta-1},
we cannot take advantage of the boundary gradient term.
So, we have to use the trace theorem as in
Lemma~\ref{moser-u}. As before, this forces us to assume some
more summability of the initial value of $\teta$.
Namely, we have
\begin{lemm}[Moser iterations for $\teta$, case $\alpha>0$, $\beta=0$]%
 \begin{sl}\label{moser-teta-2b}
 Let $S\in \RR$, $T\ge S+2$, $\alpha>0$, and $\beta=0$.
 Let $(u,\teta)$ be a couple of sufficiently smooth
 functions solving in a suitable sense system\/
 \eqref{eq-abstr}-\eqref{bc-abstr} over $(S,T)$.
 Moreover, let us assume that, for some $\epsilon\in(0,1)$,
 \eqref{tetaf2-1}-\eqref{tetaf2-2} hold, together with
 \begin{equation}\label{tetaz2b}
   \teta(S)\in L^{4+\epsilon}(\Omega),~~
          \teta(S)\in L^{4+\epsilon}(\Gamma), \qquad
    \| \teta(S) \|_{L^{4+\epsilon}(\Omega)}
     + \| \teta(S) \|_{L^{4+\epsilon}(\Gamma)} =: \Theta.
 \end{equation}
 Then, for any $\tau\in (0,\min\{1,T-S\})$,
 \eqref{st-abstr-teta} holds.
\enle
\begin{proof}
As before, here we allow $c$ to depend on $T-S$.
We take for simplicity $\alpha=1$. The analogue of
\eqref{lemma2-21} reads now
\begin{equation}\label{lemma3-21}
  \ddt \big( \| \teta \|_p^p + \| \teta \|_{p,\Gamma}^p \big)
   + \kappa \big\| \teta^{\frac{p-2}2} \big\|^2_V
     \le p \io |f| \teta^{p-1}
     + c \| \teta \|_{\frac{p-2}2}^{p-2}.
\end{equation}
Then, we take first $p=4+\epsilon$. Owing to
Lemma~\ref{moser-teta-1}, we then get
\begin{equation}\label{st-e0}
  \| \teta \|_{L^\infty(S,T;L^{4+\epsilon}(\Omega))}
   + \| \teta \|_{L^\infty(S,T;L^{4+\epsilon}(\Gamma))}
  \le Q( F, L, \Theta ).
\end{equation}
Let us now repeat \eqref{lemma3-21} for a new $p$ to be chosen below.
Let also $r=(3+\epsilon')/(2+\epsilon')$, where $\epsilon'$ is
a number, also to be chosen below, in the range $(0,\epsilon]$.
We then obtain
\begin{equation}\label{lemma3-2}
  p\io |f| \teta^{p-1}
   \le p \|f\|_{3+\epsilon'}\| \teta^{p-1} \|_r
   \le p \| f \|_{3+\epsilon} \| \teta \|_{r(p-1)}^{p-1}.
\end{equation}
Hence, integrating \eqref{lemma3-21} over $(\ttt,T)$,
for some $\ttt\ge S$ (to be chosen later) and using once
more continuity of the trace from $V$ to $L^4(\Gamma)$,
we arrive (compare with \eqref{lemma1-4}) at
\begin{align} \no
  & \| \teta \|_{L^{\infty}(\ttt,T;L^{p}(\Omega))}^p
   + \| \teta \|_{L^{\infty}(\ttt,T;L^{p}(\Gamma))}^p
   + \| \teta \|_{L^{p-2}(\ttt,T;L^{2p-4}(\Gamma))}^{p-2}
   + \| \teta \|_{L^{p-2}(\ttt,T;L^{3p-6}(\Omega))}^{p-2} \\
 \no
 & \mbox{}~~~~~~~~~~
   \le \| \teta(\ttt) \|_p^p
    + \| \teta(\ttt) \|_{p,\Gamma}^p
    + c p F \| \teta \|_{L^{2(p-1)}(\ttt,T;L^{r(p-1)}(\Omega))}^{p-1}
   + c \| \teta \|_{L^{p-2}(\ttt,T;L^{\frac{p-2}2}(\Omega))}^{p-2}\\
 \label{lemma3-4}
 & \mbox{}~~~~~~~~~~
   \le \| \teta(\ttt) \|_p^p
    + \| \teta(\ttt) \|_{p,\Gamma}^p
    + c ( p F + 1) \| \teta \|_{L^{2(p-1)}(\ttt,T;L^{r(p-1)}(\Omega))}^{p-1}.
\end{align}
%
%the third term on the \lhs\ coming from continuity of
%the trace operator from $V$ to $L^4(\Gamma)$.
Now, we perform Moser iterations as in
the proof of Lemma~\ref{moser-u}, starting from
$p_0=4+\epsilon$. Then, in place of \eqref{interpol12}, we
get the system
\begin{equation} \label{interpol22}
  \displaystyle
   \begin{cases}
     \frac{1-\rho}{p_i-2} = \frac{1}{2(p_{i+1}-1)},\\
      \frac{\rho}{p_i} + \frac{1-\rho}{3(p_i-2)}\,=\,\frac{1}{r(p_{i+1}-1)},
   \end{cases}
\end{equation}
whence one computes, similarly as before,
\begin{equation} \label{interpol25}
  K_{\epsilon'}=\frac{9+5\epsilon'}{9+3\epsilon'}.
\end{equation}
and, finally,
\begin{align} \no
  p_{i+1} & = \frac12 \frac{p_i-2}{1-\rho} + 1
    = \frac12 \Big(1 + \frac{p_i}{p_i-2} K_{\epsilon'} \Big) (p_i-2) + 1\\
  \label{interpol28}
  & = \frac{ K_{\epsilon'} + 1 }2 p_i = \frac{9+4\epsilon'}{9+3\epsilon'} p_i =: H p_i
  = H^{i+1} p_0.
\end{align}
Note that, as before, $H$ is larger than $1$ since $\epsilon'>0$.

Then, we can follow with minor variations
the proof of Lemma~\ref{moser-u}
up to formula \eqref{moser10}.
In particular, $J_i$ is now defined by setting
\begin{equation}\label{Jp-new}
  J_i^{p_i} := \| \teta \|_{L^{\infty}(\tau_i,T;L^{p_i}(\Omega))}^{p_i}
   + \| \teta \|_{L^{\infty}(\tau_i,T;L^{p_i}(\Gamma))}^{p_i}
   + \| \teta \|_{L^{p_i-2}(\tau_i,T;L^{2p_i-4}(\Gamma))}^{p_i-2}
   + \| \teta \|_{L^{p_i-2}(\tau_i,T;L^{3p_i-6}(\Omega))}^{p_i-2}.
\end{equation}
Then, the analogue of \eqref{moser10b} reads
\begin{align} \no
 & \| \teta(\tau_{i+1}) \|_{p_{i+1}}^{p_i-2}
  + \| \teta(\tau_{i+1}) \|_{p_{i+1},\Gamma}^{p_i-2}
  \le \frac1{\ttt_{i+1}}\int_{\tau_i}^{\tau_i+\ttt_{i+1}}
   \Big( \| \teta(t) \|_{p_{i+1}}^{p_{i}-2}
    + \| \teta(t) \|_{p_{i+1},\Gamma}^{p_{i}-2} \Big)\,\dit\\
  \label{moser10b-new}
  & \mbox{}~~~~~~~~~~
  \le c \frac{i^2}{\tau} \int_{\tau_i}^{\tau_i+\ttt_{i+1}}
    \Big( \| \teta(t) \|_{3(p_i-2)}^{p_{i}-2}
    + \| \teta(t) \|_{2(p_i-2),\Gamma}^{p_i-2} \Big)\,\dit
  \le c \frac{i^2}{\tau} J_i^{p_i}. 
\end{align}
Hence, to control the \rhs\ of \eqref{lemma3-4}, we also need to
take the $\frac{p_i}{p_i-2}$-th power of the above inequality. It
is readily checked that this does not affect 
the validity of the argument since the exponent
$\frac{p_i}{p_i-2}$ is asymptotically close to $1$. Note also that
\eqref{moser10b-new} holds provided that $p_{i+1}\le 2(p_i-2)$.
As $p_{i+1}$ is given by \eqref{interpol28},
this means that we need
\begin{equation} \label{condsu-p}
  (4+\epsilon) H^{i+1}
   \le 2 \big( (4+\epsilon) H^i - 2 \big).
\end{equation}
which is easily shown to be true for every $i\ge 0$ provided
that we choose $\epsilon'\in(0,\epsilon]$ so small that
\begin{equation} \label{condsu-p-2}
  (4+\epsilon) H
   = (4+\epsilon) \frac{9+4\epsilon'}{9+3\epsilon'}
   \le 2 \big( (4+\epsilon) - 2 \big).
\end{equation}
Actually, also the choice $\epsilon'=\epsilon$ works.
Notice that, since we used $p_i < p_{i+1}\le 2(p_i-2)$,
it follows that the regularity \eqref{tetaz2b} is
needed in order for the iteration scheme to work.
From this point on, the proof proceeds once again similarly
with that of Lemma~\ref{moser-u}, hence we can omit the details.
\end{proof}
\noindent%
As before, we have a better result in case the ``initial''
value of $\teta$ is uniformly bounded. The proof follows the
preceding ones up to straightforward modifications.
\beco\label{moser-teta-3}
 Let $S\in \RR$, $T\ge S+2$, $\alpha\ge 0$, $\beta\ge 0$.
 Let $(u,\teta)$ be a couple of sufficiently smooth
 functions solving in a suitable sense system\/
 \eqref{eq-abstr}-\eqref{bc-abstr} over $(S,T)$.
 Let also \eqref{tetaf2-1}-\eqref{tetaf2-2} hold,
 together with
 \begin{equation}\label{tetaz3}
   \teta(S)\in L^{\infty}(\Omega),~~
         \alpha \teta(S)\in L^{\infty}(\Gamma), \qquad
    \| \teta(S) \|_{L^{\infty}(\Omega)}
     + \alpha \| \teta(S) \|_{L^{\infty}(\Gamma)} =: \Theta.
 \end{equation}
 Then, we have 
 \begin{equation}\label{st-teta-unif}
   \| \teta \|_{L^\infty((S,T)\times \Omega)} \le
    Q\big( L, F, \Theta, T-S \big).
 \end{equation}
\enco
\beos 
 The $L^{3+\epsilon}$ space regularity of $f$,
 stated in \eqref{hypfw} with $\epsilon>0$, is strongly exploited 
 in the proofs (cf.~\eqref{interpol18} and \eqref{interpol25}). 
 Indeed, if one weakens it by putting $\epsilon=0$, 
 no additional summability is gained
 at subsequent iteration steps and the argument fails.
\eddos

%%%%%%%%%%%%%%%%%%%%%%%%%%%%%%%%%%%%%%%%%%%%%%%%%%%%%%%%%%%%%%%%%%%

\subsection{Further regularity of time derivatives}
\label{subsec:regtime}

We now prove that, for smooth solutions of Problem~(P), the
boundedness of $\teta$ and $u$ implies some $L^p$-regularity
of time derivatives at least for strictly positive times.
The proof of this result is a little bit tricky in
the case $\beta>0$ since the diffusion operators on
$\Omega$ and on $\Gamma$ act on different functions.
In particular, we cannot treat the case when $\beta>0$
and $\alpha=0$ because in the derivation of the
estimates the contribution of the boundary term $\teta_t$
is explicitly needed.
\bele\label{lemma-time}
 Let $S\in\RR$, $T>S$, $\alpha\ge 0$, $\beta\ge 0$.
 % and let\/ \eqref{inpiuf} hold for some $\epsilon\in(0,1)$.
 If $\beta>0$, then let $\alpha>0$.
 Let $(u,\teta)$ be a couple of sufficiently
 smooth functions solving in a suitable sense
 system\/ \eqref{eq-abstr}-\eqref{bc-abstr}
 over the time interval $(S,T)$.
 Moreover, %let us assume\/
 %\eqref{eq-abstr}-\eqref{bc-abstr} and
 let
 \begin{align}\label{l4-41new1}
   & u \in L^2(S,T;V),~~\beta \teta \in L^2(S,T;V\Ga), \qquad
    \| u \|_{L^2(S,T;V)} + \beta \| \teta \|_{L^2(S,T;V\Ga)} \le U,\\
  \label{l4-41new2}
   & f \in L^2(S,T;H), \qquad
    \| f \|_{L^2(S,T;H)} =: F,\\
  \label{l4-42new}
   & u,\teta \in L^\infty((S,T)\times \Omega), \qquad
     \| u \|_{L^\infty((S,T)\times \Omega)}
     + \| \teta \|_{L^\infty((S,T)\times \Omega)} \le U,
  \end{align}
 for some (given) constants $F>0$, $U>0$.
 Then, for any $\tau\in(0,\min\{1,T-S\})$, we have 
 \begin{equation}\label{st-tetatnew}
   \| \teta_t \|_{L^2(S+\tau,T;H)}
    + \alpha \| \teta_t \|_{L^2(S+\tau,T;H\Ga)}
    + \| u \|_{L^\infty(S+\tau,T;V)}
    + \beta \| \teta \|_{L^\infty(S+\tau,T;V\Ga)}
   \le Q\big(F,U,\tau^{-1},T-S\big). 
 \end{equation}
\enle
\begin{proof}
We just consider the case when $\alpha,\beta>0$, which is
the most difficult one. For simplicity, we can then
set $\alpha=\beta=1$ and write
\begin{equation}\label{ibaro}
  \ibaro v:= \io v + \iga v,
\end{equation}
whenever $v$ is, say, an element of $\calH$. Then, we
test both \eqref{eq-abstr} and \eqref{bc-abstr} by $u_t=\teta_t/\teta^2$
(note that we are always assuming the solution to be smooth
enough for our purposes). Then, standard integrations by
parts lead to
\begin{equation}\label{lemma4-11}
  \frac12 \ddt \Big( \| \nabla u \|^2
   + \iga \frac{|\nabla\Ga \teta|^2}{\teta^2} \Big)
   + \ibaro \frac{\teta_t^2}{\teta^2}
  \le \iga \frac{\teta_t|\nabla\Ga \teta|^2}{\teta^3}
  + \io f \frac{\teta_t}{\teta^2}.
\end{equation}
Of course, using \eqref{l4-42new},
the latter term can be simply estimated as follows:
\begin{equation}\label{lemma4-12}
  \io f \frac{\teta_t}{\teta^2}
   \le \frac14\io\frac{\teta_t^2}{\teta^2}
   + Q(U) \| f \|^2.
\end{equation}
The first term on the \rhs\ of \eqref{lemma4-11} is
more delicate. Actually, owing to standard interpolation
inequalities, we have
\begin{align} \no
  \iga \frac{\teta_t|\nabla\Ga \teta|^2}{\teta^3}
   & \le Q(U) \| \teta_t \|\Ga \| \nabla\Ga \teta \|_{4,\Gamma}^2
   \le Q(U) \| \teta_t \|\Ga \| \teta \|_{H^2(\Gamma)} \| \teta \|_{\infty,\Gamma} \\
  \no
   & \le Q(U) \big( \| \teta_t \|\Ga^2 + \| - \Delta\Ga \teta \|\Ga^2
    + \| \teta \|\Ga^2 \big) \\
  \label{lemma4-13}
   & \le Q(U) \big( \| \teta_t - \Delta\Ga \teta \|\Ga^2 \big)
    + Q(U) - Q(U) \ddt \| \nabla\Ga \teta\|\Ga^2.
\end{align}
Then, the last term is moved to the \lhs\ and will give
a positive contribution. We now estimate the first term on
the \rhs, which requires some work.
Actually, we first notice that, comparing terms
in equation~\eqref{bc-abstr}, using
the trace theorem and elliptic regularity results
(cf., e.g., \cite[Thm.~2.5.13 and Thm.~2.7.4]{BrG}), 
for arbitrarily small 
but otherwise fixed $\delta\in(0,1/2)$,
we obtain
\begin{align} \no
  Q(U) \big( \| \teta_t - \Delta\Ga \teta \|\Ga^2 \big)
  +  Q(U) \| \dn u \|\Ga^2
  & \le Q(U) \| \dn u \|\Ga^2
   \le Q(U) \| u \|_{H^{\frac32+\delta}(\Omega)}^2 \\
 \no
  & \le Q(U) \Big( \| u \|^2
    + \| \Delta u \|_{H^{-\frac12+\delta}(\Omega)}^2
    + \| u \|_{H^{1+\delta}(\Gamma)}^2 \Big)  \\
 \label{lemma4-14}
  & \le Q(U) \Big( 1
    + \| \Delta u \|^p \| u \|^{2-p}
    + \| u \|_{H^{1}(\Gamma)}^{2-2\delta} \| u \|_{H^{2}(\Gamma)}^{2\delta} \Big),
\end{align}
for a suitable $p\in(0,2)$ depending on the
choice of $\delta$. Then, using~\eqref{eq-abstr}
and Young's inequality, we obtain
\begin{align} \no
  Q(U) \| \Delta u \|^p \| u \|^{2-p}
   & \le Q(U) \| \Delta u \|^p
   \le Q(U) \| \teta_t - f \|^p \\
 \label{lemma4-15}
  & \le \sigma \io \frac{\teta_t^2}{\teta^2}
   + Q_\sigma(U) \big( 1 + \| f \|^2 \big).
\end{align}
Moreover, comparing once more terms in \eqref{bc-abstr}
and noting in particular that
\begin{equation}\label{lemma4-15b}
  \Delta\Ga u = \dive\Ga \big( u^2 \nabla\Ga \teta \big),
\end{equation}
we arrive at
\begin{align} \no
  \| \Delta\Ga u \|\Ga & \le Q(U) \| \Delta\Ga \teta \|\Ga
   +  Q(U) \| \nabla\Ga \teta \|_{4,\Gamma}^2\\
 \no
  & \le Q(U) \| \Delta\Ga \teta \| + Q(U)
   + Q(U) \| \Delta\Ga \teta \|_{\Gamma} \| \teta \|_{\infty,\Gamma} \\
 \label{lemma4-15c}
  & \le Q(U) \| \Delta\Ga \teta \|\Ga + Q(U).
\end{align}
Hence, the last term in \eqref{lemma4-14} can be
controlled this way:
\begin{align} \no
  Q(U) \| u \|_{H^{1}(\Gamma)}^{2-2\delta} \| u \|_{H^{2}(\Gamma)}^{2\delta}
    & \le Q_\sigma(U) \| u \|_{H^{1}(\Gamma)}^2
     + \sigma \| \Delta\Ga u \|\Ga^2 \\
 \no
  & \le Q_\sigma(U) \| u \|_{H^{1}(\Gamma)}^2
  + \sigma Q(U) \| \Delta\Ga \teta \|\Ga^2
  + \sigma Q(U) \\
 \no
  & \le  Q_\sigma(U) \| u \|_{H^{1}(\Gamma)}^2
   + \sigma Q(U) \| \teta_t \|\Ga^2
   + \sigma Q(U) \| \dn u \|\Ga^2
   + \sigma Q(U) \\
 \label{lemma4-16}
  & \le Q_\sigma(U) \| \teta \|_{H^{1}(\Gamma)}^2
   + \sigma Q(U) \iga \frac{\teta_t^2}{\teta^2}
   + \sigma Q(U) \| \dn u \|\Ga^2 + \sigma Q(U),
\end{align}
where $\sigma>0$ is small and $Q_\sigma$ also
depends on $\sigma$. In particular, taking
$\sigma$ small enough both in \eqref{lemma4-15}
and in \eqref{lemma4-16}, and collecting
\eqref{lemma4-12}-\eqref{lemma4-16},
\eqref{lemma4-11} gives
\begin{align} \no
  & \frac12 \ddt \Big( \| \nabla u \|^2
   + \iga \frac{|\nabla\Ga \teta|^2}{\teta^2}
   + Q(U) \| \nabla\Ga \teta\|\Ga^2 \Big)
   + \frac12 \ibaro \frac{\teta_t^2}{\teta^2} \\
 \label{lemma4-17}
  & \mbox{}~~~~~
   \le Q(U) \big( \| f \|^2
   +  \| \teta \|_{H^{1}(\Gamma)}^2
   + 1 \big).
\end{align}
Now, thanks to \eqref{l4-41new1}, for any $\tau$
as in the statement, we
can choose $\ttt\in (S,S+\tau)$ such that
\begin{equation}\label{lemma4-18}
  \| \nabla u (\ttt) \|^2
   + \| \nabla\Ga \teta (\ttt) \|\Ga^2 \le Q(U,\tau^{-1}).
\end{equation}
The thesis then follows by integrating \eqref{lemma4-17}
over $(\ttt,T)$.
\end{proof}
\noindent%
As before, for regular ``initial'' data, the property
holds starting from the initial time:
\beco\label{lemma-time-2}
 Let $S\in\RR$, $T>S$, $\alpha\ge 0$, $\beta\ge 0$.
 % and let\/ \eqref{inpiuf} hold for some $\epsilon\in(0,1)$.
 If $\beta>0$, then let $\alpha>0$.
 Let $(u,\teta)$ be a couple of sufficiently
 smooth functions solving in a suitable sense
 system\/ \eqref{eq-abstr}-\eqref{bc-abstr}
 over the time interval $(S,T)$.
 Moreover, let\/ \eqref{l4-41new1}-\eqref{l4-42new}
 hold together with
 \begin{equation}\label{l4-51}
    u(S) \in V, \quad \beta \teta(S) \in V\Ga,
     \qquad \| u(S)\|_V + \beta \| \teta(S) \|_{V\Ga} =: Z,
  \end{equation}
 for some $Z>0$. Then,
 \begin{equation}\label{st-tetatnew2}
   \| \teta_t \|_{L^2(S,T;H)}
    + \alpha \| \teta_t \|_{L^2(S,T;H\Ga)}
    + \| u \|_{L^\infty(S,T;V)}
    + \beta \| \teta \|_{L^\infty(S,T;V\Ga)}
   \le Q\big(F,U,Z,T-S). 
 \end{equation}
\enco

%%%%%%%%%%5%%%%%%%%%%%%%%%%%%%%%%%%%%%%%%%%%%%%%%%%%%%%%%%%%%%%%%%%%%%%%%%%%%%

\section{Existence for smooth data}
\label{sec:smooth}

Based on the previous lemmas, we will prove here that,
for sufficiently regular initial data and source term, Problem~(P)
has a global solution in a rather good regularity class.

%%%%%%%%%%5%%%%%%%%%%%%%%%%%%%%%%%%%%%%%%%%%%%%%%%%%%%%%%%%%%%%%%%%%%%%%%%%%%%

\subsection{Local existence of smooth solutions}
\label{sec:loc:smooth}

We start by showing that, for regular data, a local
in time smooth solution exists. For the sake of
simplicity, we give the proof only in the case when
both $\alpha$ and $\beta$ are strictly positive. Actually,
the other cases can be treated with differences
that are mainly of technical character.
\bete\label{teo:loc:reg}
 Let $\alpha>0$, $\beta>0$. Let $T>0$ %, $a\in(0,1)$
 and let
 \begin{align} \label{hpf:reg}
   & f \in C^{0}([0,T]\times \barO),\\
  \label{hpteta:reg}
   & \teta_0 \in H^2(\Omega),
    \qquad \tetagiu\le \teta_0(x)\le \tetasu~~\text{for all }\,x\in\barO,
%  \label{hpu:reg}
%   u_0:=-\frac1{\teta_0} \in C^{0,a}(\barO),
%    \qquad |\u_0(x)| \ge \ugiu> 0~~\text{for all }\,x\in\barO.
 \end{align}
 for suitable constants $0<\tetagiu<\tetasu$. Then, there exists
 $T_0\in(0,T]$ depending on $f$ and $\teta_0$ (and in particular
 on $\tetagiu,\tetasu$), such that\/ {\rm Problem~(P)}
 admits a solution $(\teta,u)$ over the time interval $[0,T_0]$
 satisfying
 \begin{align} \label{reg:teta:reg}
   & \teta\in H^1(0,T_0;\calH) \cap C^0([0,T_0]\times \barO),
    \qquad \frac{\tetagiu}2\le \teta(t,x)\le 2\tetasu~~\text{for all }\,(t,x)\in [0,T_0]\times \barO,\\
  \label{reg:u:reg}
   & u \in L^2(0,T_0;H^2(\Omega)),~~\teta\in L^2(0,T_0;H^2(\Gamma)).
 \end{align}
\ente
\begin{proof}
 We only sketch it, since it follows from more or less classical
 arguments for quasilinear parabolic problems. The key idea is to
 regularize the monotone function $\gamma(r)=-1/r$ by introducing
 a function $\gamma_R:\RR\to \RR$ with the following properties
 \begin{align} \label{gammaR1}
   \gamma_R \in C^2(\RR), \quad \gamma_R''\in L^\infty(\RR),
   \qquad \gamma_R(r) = \gamma(r)~~\perogni r \in \Big[\frac{\tetagiu}2,2\tetasu\Big],\\
  \label{gammaR2}
    \esiste 0<m_R<M_R\quext{such that }\, m_R \le \gamma_R'(r)\le M_R~~%
     \perogni r\in\RR.
 \end{align}
 In other words, $\gamma_R$ coincides with $\gamma$ over
 an interval that is strictly larger than the range of the initial
 values of the problem. Outside that interval, we substitute
 it by a smooth and bi-Lipschitz approximation
 (e.g., we can take the first order Taylor expansion
 and then mollify).
 At this point, we can solve the problem by means of a fixed point
 argument. We just give the highlight of
 this procedure and set, for simplicity, $\alpha=\beta=1$.

 As a first step, we take a function $\phi$
 in the closed ball $\Phi$ of radius $1$ of the space
 \begin{equation} \label{newspace}
   L^2(0,T_0;H^{\epsi}(\Gamma)),
 \end{equation}
 where $\epsi\in(0,1/4)$ 
 is fixed but otherwise arbitrary.
 Then, we consider the initial-value problem
 \begin{equation} \label{probound}
   \eta_t - \Delta\Ga \eta = \phi, \quad
     \eta|_{t=0} = \teta_0|\Ga,
     \qquext{on }\,\Gamma.
 \end{equation}
 and we claim that it admits a unique solution $\eta$ lying
 in a bounded closed ball of the space
 \begin{equation} \label{regbound}
    H^1(0,T_0;H^{\epsi}(\Gamma)) \cap 
     L^2(0,T_0;H^{2+\epsi}(\Gamma)).
 \end{equation}
 The simpler way to obtain this property is probably to 
 prove it directly by interpolation. Actually, the regularity
 \eqref{newspace} of the datum $\phi$ implies the regularity
 \eqref{regbound} of the solution $\eta$ for $\epsi=0,1$ 
 (this is classical and can be proved easily; notice that
 $\Gamma$ is a smooth manifold without boundary). Moreover,
 the solution operator $\phi\mapsto \eta$ is continuous 
 for $\epsi=0,1$. Hence, by standard interpolation results
 we get continuity of the solution map for $\epsi\in(0,1/4)$.
 Notice also that, for such a choice of $\epsi$, the regularity
 of the boundary initial datum (which is known to lie in $H^{3/2}(\Gamma)$
 in view of \eqref{hpteta:reg} and the trace theorem), is sufficient.
 Next, we observe that, again by interpolation, the space
 \eqref{regbound} is continuously embedded into $C^{0,a}([0,T_0]\times \Gamma)$
 for some $a>0$ depending on the choice of $\epsi$
 (recall that $\Omega$ is smooth and its boundary
 $\Gamma$ is a 2-dimensional manifold). Note also that this property
 may fail for $\epsi=0$. 
 Thus, we can take $T_0$ so small that, additionally,
 \begin{equation} \label{limtrace}
   \qquad \frac{3\tetagiu}4\le \eta(t,x)\le \frac32\tetasu
     \quext{for all }\,(t,x)\in [0,T_0]\times\Gamma.
 \end{equation}
 As a subsequent step, we solve the problem
 \begin{align} \label{proomega1}
   & \teta_t - \Delta \gamma_R(\teta) = f, \quad
     \teta|_{t=0} = \teta_0, \quext{in }\,\Omega,\\
  \label{proomega2}
   & \gamma_R(\teta) = \gamma_R(\eta), \quext{on }\,\Gamma.
 \end{align}
 This is a quasilinear parabolic system with bi-Lipschitz
 continuous nonlinearity $\gamma_R$ and Dirichlet boundary conditions.
 Then, by standard methods one may verify that it admits
 a unique solution $\teta$ such that
 \begin{equation} \label{regomega1}
   \teta\in H^1(0,T_0;H), \quad
    \gamma_R(\teta)\in L^\infty(0,T_0;H^1(\Omega)) \cap 
     L^2(0,T_0;H^2(\Omega)).
 \end{equation}
 For instance, the basic a priori estimate corresponding
 to the above regularity can be obtained testing
 the equation in \eqref{proomega1} by the function
 $\partial_t(\gamma_R(\teta) - \calR(\gamma_R(\eta)))$, where
 $\calR$ denotes the harmonic extension operator, namely
 \begin{equation} \label{deficalR}
   -\Delta \calR(v) = 0~~\text{in }\,\Omega, \qquad
    \calR(v) = v~~\text{on }\,\Gamma,
 \end{equation}
 where $v$ is, say, a function in $H\Ga$. Note that
 the regularity \eqref{regbound} of the trace
 and the $C^2$-regularity of $\gamma_R$ are essential for the
 sake of obtaining \eqref{regomega1}, as a direct check
 permits us to verify.

 Then, being $\gamma_R$ bi-Lipschitz, 
 $\partial_t(\gamma_R(\teta))$ has the same $L^2$-regularity 
 as $\teta_t$ (cf.~\eqref{regomega1}). Hence, we may conclude that
 $\gamma_R(\teta)$ lies
 in a bounded closed ball of the space
 \begin{equation} \label{regomeg}
   H^1(0,T_0;H) \cap L^\infty(0,T_0;H^1(\Omega)) \cap L^2(0,T_0;H^2(\Omega)).
 \end{equation}
 Moreover, the equation in \eqref{proomega1} is
 quasilinear and uniformly parabolic and
 it has uniformly bounded initial data
 (by \eqref{hpteta:reg}), boundary data
 (by \eqref{limtrace}), and forcing term
 (by \eqref{hpf:reg}). Thus, standard barrier
 arguments entail that, up to possibly taking
 a smaller initial time $T_0$ (in a way that only depends
 on the known norms of $\teta_0$ and $f$, on the
 regularized function $\gamma_R$, and on
 the truncation values $\tetagiu,\tetasu$),
 \begin{equation} \label{regomeg2}
   \frac{\tetagiu}2\le \teta(t,x)\le 2\tetasu~~
     \text{for all }\,(t,x)\in [0,T_0]\times\barO.
 \end{equation}
 Next, let us notice that, by interpolation of Sobolev spaces
 (cf.~once more \cite[Chap.~2]{BrG} for details),
 the space in \eqref{regomeg} is continuously embedded into
 $H^\delta(0,T_0;H^{2-2\delta}(\Omega))$ for any $\delta\in(0,1)$.
 Hence, at least for $\delta\in(0,1/4)$, we may apply
 the trace theorem \cite[Thm.~2.7.4]{BrG}. In particular,
 taking $\delta=1/8$, we obtain that the function $\dn\gamma_R(\teta)$ 
 lies in a bounded closed ball
 of the space
 \begin{equation} \label{regbound3}
    H^{1/8}(0,T_0;H^{1/4}(\Omega)),
 \end{equation}
 which, in view of the fact that $\epsi\in(0,1/4)$, is continuously and
 compactly embedded into the space in \eqref{newspace}.
 Hence, the map
 \begin{equation} \label{fixpo}
   \calT: \Phi \to L^2(0,T_0;H^{\epsi}(\Gamma)), \qquad
   \calT: \phi \mapsto \dn\gamma(\teta),
 \end{equation}
 is compact. Next, we prove that, for $T_0$ small enough,
 $\calT$ takes values into $\Phi$. Actually, interpolating the latter
 two spaces in \eqref{regomeg}, we obtain that
 $\gamma_R(\teta)$ lies in a bounded closed ball
 of $L^{16/7}(0,T_0;H^{7/4}(\Omega))$, whence,
 applying once more \cite[Thm.~2.7.4]{BrG},
 $\dn\gamma_R(\teta)$ also lies in a bounded closed ball
 of the space $L^{16/7}(0,T_0;H^{1/4}(\Gamma))$. Moreover, the
 radius of that ball has the form $Q_1(\|\teta_0\|_{H^2(\Omega)})$.
 In particular, the radius is independent of $T_0$ in view of the fact
 that all the constants appearing in the estimates performed so far 
 are also independent of $T_0$ since they come from parabolic 
 regularity theorems. In other words, we do not need to use
 Sobolev's embeddings (in which case the embedding constants may explode
 as the measure of the time interval $(0,T_0)$ becomes small).
 Hence, since $16/7>2$ and $\epsi\in (0,1/4)$, comparing with 
 \eqref{newspace}, we can take $T_0$ small enough so that
 $\calT$ takes values into $\Phi$. Finally, by a number of standard checks
 one may verify that $\calT$ is continuous. Hence, the Schauder fixed point
 theorem may be applied to $\calT$. This provides a local solution
 to the system
 \begin{align}\label{eq-rego}
   & \teta_t - \Delta \gamma_R(\teta) = f,
    \quext{in }\,\Omega,\\
  \label{bc-rego}
   & \eta_t - \Delta\Ga \eta = - \dn \gamma_R(\teta),
   \quext{on }\,\Gamma,
 \end{align}
 plus the initial conditions and the boundary condition
 $\teta|\Ga=\eta$. However, thanks to \eqref{regomeg2}
 and the latter \eqref{gammaR1}, $\gamma_R(\teta)$ coincides
 with $\gamma(\teta)=-1/\teta=:u$ everywhere in
 $[0,T_0]\times \barO$. Hence, $(\teta,u)$ is a
 local smooth solution to~\eqref{eq-abstr}-\eqref{bc-abstr}.
 This concludes the proof.
\end{proof}
%

%%%%%%%%%%5%%%%%%%%%%%%%%%%%%%%%%%%%%%%%%%%%%%%%%%%%%%%%%%%%%%%%%%%%%%%%%%%%%%

\subsection{Energy estimate}
\label{sec:esti}

In this section we prove the basic
energy estimate satisfied by solutions
of Problem~(P).
We start by recalling a generalized
version of Poincar\'e's inequality
(see for instance \cite[Lemma~3.2]{SSZ1} for a proof).
\bele\label{lemma-log-poincare}
 Assume $\Omega$ is a bounded open subset of $\mathbb{R}^d$. Suppose
 $v\in W^{1,1}(\Omega)$ and $v\ge 0$ a.e. in $\Omega$.
 Then, setting $K:=\int_{\Omega}(\log v)^{+}$,
 the following estimate holds:
 \begin{equation}\label{log-poincare}
   \|v\|_{L^1(\Omega)} \le |\Omega| e ^{C_1 K}
    + \frac{C_2}{|\Omega|}\|\nabla v\|_{L^1(\Omega)},
 \end{equation}
 the positive constants $C_1$ and $C_2$ depending only on $\Omega$.
\enle
\noindent
\begin{lemm}[Energy estimate]\begin{sl}\label{lemma:energy}
 Let $\alpha\ge 0$, $\beta\ge 0$ and let\/
 \eqref{hpteta0w}-\eqref{hpeta0w} and\/ \eqref{hypfw} hold.
 Let $(\teta,u)$ be a sufficiently smooth solution to\/
 {\rm Problem~(P)}. Then, we have 
 \begin{equation}\label{stimau}
   \calE(t) +
    \int_0^t \big( \| \nabla u \|^2
     + \beta \| \nabla\Ga \log \teta \|^2\Ga \big)
     + \| u \|_{L^2(t,t+1;V)}^2 \le 
     Q\big(\EE_0, \| f \|_{L^2(0,t;L^{6/5}(\Omega))} \big), 
    \quad\perogni t\ge 0.
 \end{equation}
\enle
\begin{proof}
Test \eqref{calore} and \eqref{calorega} with $1 + u$. This gives
\begin{equation}\label{energyprovv}
  \ddt \calE(t) +
   \big( \|\nabla u\|^2 + \beta \|\nabla\Ga \log \teta\|\Ga^2 \big)
   = \io f (1+u).
\end{equation}
Then, using \eqref{hypfw} and the Poincar\'e-Wirtinger
inequality, we estimate the \rhs\ as follows:
\begin{equation}\label{energyprovv2}
  \io f (1+u) = \io f u
   = \io f \big( u - m\OO(u) \big)
   \le \frac12 \| \nabla u \|^2 + c \| f \|_{L^{6/5}(\Omega)}^2,
\end{equation}
with $c$ independent of time. 
Integrating \eqref{energyprovv} in time from $0$ to $t$, we obtain
\begin{equation}\label{energyeq}
  \calE(t) + \int_0^t
   \Big( \frac12\|\nabla u\|^2 + \beta \|\nabla\Ga \log \teta\|\Ga^2 \Big)
   \le \EE_0 + c \| f \|_{L^2(0,t;L^{6/5}(\Omega))}^2,
    \quad\perogni t\ge 0.
\end{equation}
To obtain the control on the full $V$-norm of $u$, we use
Lemma~\ref{lemma-log-poincare}. Thus, taking $v=1/\teta=-u$ in
\eqref{log-poincare}, we infer
\begin{equation}\label{poinc-u}
  \|u\|_{L^1(\Omega)} \le \disp |\Omega|e^{C_1 \io\log^-\teta}
   + \frac{C_2}{|\Omega|}\|\nabla u\|_{L^1(\Omega)}.
\end{equation}
Hence, squaring, integrating over $(t,t+1)$, and using \eqref{energyeq},
we obtain
\begin{equation}\label{st21}
  \| u \|_{L^2(t,t+1;L^1(\Omega))} \le Q(\EE_0).
\end{equation}
The thesis follows combining \eqref{st21} and
\eqref{energyeq}.
\end{proof}
%

%%%%%%%%%%5%%%%%%%%%%%%%%%%%%%%%%%%%%%%%%%%%%%%%%%%%%%%%%%%%%%%%%%%%%%%%%%%%%%

\subsection{Proof of global existence}
\label{sec:global}

Our next aim is to prove that, for smooth initial
data, the solution constructed in Section~\ref{sec:loc:smooth}
has, in fact, a global in time character. We can treat
all cases with the exception of $\alpha=0$, $\beta>0$.
This case will be dealt with separately in
Section~\ref{subsec:asymteta} below.
\bete\label{teo:glob:reg}
 Let $\alpha\ge 0$, $\beta\ge 0$ and, if $\alpha=0$, then
 let $\beta=0$. Let \eqref{hypfw} hold
 together with
 \begin{equation}\label{hp:uteta:reg}
   \teta_0 \in \calV,
      \qquad \tetagiu\le \teta_0(x)\le \tetasu~~\text{for all }\,x\in\barO,
 \end{equation}
 for some $0<\tetagiu<\tetasu$.
 Then, there exists a global solution to\/
 {\rm Problem (P)} satisfying, for all $T>0$,
 \begin{align} \label{reg:teta:reg2}
   & \teta\in H^1(0,T;H) \cap L^\infty((0,T)\times\Omega)
    \cap L^\infty(0,T;V), \\
  \label{reg:u:reg2}
   & u \in H^1(0,T;H) \cap L^\infty((0,T)\times\Omega)
    \cap L^\infty(0,T;V) \cap L^2(0,T;H^{3/2}(\Omega)), \\
  \label{reg:eta:reg2}
   & \alpha \eta \in H^1(0,T;H\Ga), \qquad \beta \eta \in L^\infty(0,T;V\Ga)
    \cap L^2(0,T;H^2(\Gamma)).
 \end{align}
 Moreover, if either $\alpha>0$ and $\beta>0$ or $\alpha=\beta=0$,
 then we have more precisely
 \begin{equation}\label{reg:teta:reg3}
   \teta \in L^2(0,T;H^{2}(\Omega)), \quad
    u \in L^2(0,T;H^{2}(\Omega)).
  \end{equation}
\ente
\begin{proof}
Let us fix (an arbitrary) $T>0$. In this proof we will
allow the constants $c$ to depend on $T$. Indeed, we are working
on a fixed interval and are not looking for uniform estimates
at this level. 
In order to get a
local solution via Theorem~\ref{teo:loc:reg} we need to
construct sequences of regularized initial and source data
$\{\teta\zzn\}$ and $\{f_n\}$ satisfying, for all
$n\in \NN$, \eqref{hpf:reg} and, respectively,
\eqref{hpteta:reg}. Moreover, we need $\teta\znn\to\teta_0$
and $f_n\to f$ in proper ways. The details of
the regularization procedure are sketched
in Section~\ref{sec:approx} below, to which we refer
the reader.

Then, thanks to Theorem~\ref{teo:loc:reg},
for all $n>0$ we have a solution to the $n$-Problem
defined on some interval $(0,T\zzn)$,
with $T\zzn\le T$. We now deduce global in time estimates
independent of $n$ and, for the sake of simplicity,
we shall directly work on the time interval $(0,T)$
rather than on $(0,T\zzn)$. As usual, this can be justified
{\sl a posteriori}\/ by means of standard extension arguments.
To be more precise, one may first extend the approximate
solution up to the final time $T$ by applying the uniform estimates
at {\sl fixed} $n$. Indeed, these estimates imply that, 
if $T_{0,n}'\ge T_{0,n}$ is the maximal existence time of
the approximate solution in the regularity class of 
Theorem~\ref{teo:loc:reg}, then the approximate 
solution cannot explode as $t\nearrow T_{0,n}'$ (because the 
estimates are uniform) and consequently $T_{0,n}'$ must
coincide with the reference time $T$. Once the approximate
solutions are proved to exist over $(0,T)$, then one may
apply once more the estimates in order to let $n\nearrow\infty$
and remove the approximation. 

That said, we can first apply Lemma~\ref{lemma:energy},
which gives the bound
\begin{equation} \label{sti11}
  \| \teta_n\|_{L^\infty(0,T;L^1(\Omega))}
   + \alpha \| \eta_n\|_{L^\infty(0,T;L^1(\Gamma))}
   + \| u_n \|_{L^2(0,T;V)}
   + \beta \| \log \eta_n \|_{L^2(0,T;V\Ga)}
    \le c,
\end{equation}
for $c>0$ independent of $n$. Indeed, from
this point on, we go back to the notation $\eta_n$ when
we indicate the boundary value of $\teta_n$.
Next, we can apply
Corollaries~\ref{moser-u-2}, \ref{moser-teta-3}
(with $S=0$), which give
\begin{equation} \label{sti12}
  \| u_n \|_{L^\infty((0,T)\times\Omega)}
   + \| \teta_n \|_{L^\infty((0,T)\times\Omega)}
    \le c.
\end{equation}
Clearly, the same uniform boundedness properties
hold also for the traces on $\Gamma$.

Next, by Corollary~\ref{lemma-time-2}
(here the restriction on $\alpha$ and $\beta$
comes into play), we have
\begin{equation}\label{sti13}
  \| (\teta_n)_t \|_{L^2(0,T;H)}
   + \alpha \| (\eta_n)_t \|_{L^2(0,T;H\Ga)} 
   + \| u_n \|_{L^\infty(0,T;V)}
   + \beta \| \eta_n \|_{L^\infty(0,T;V\Ga)}
   \le c.
\end{equation}
Comparing terms in equation \eqref{eq-abstr},
we also obtain
\begin{equation}\label{sti14}
  \| \Delta u_n \|_{L^2(0,T;H)}
   \le c.
\end{equation}
Moreover, in the case when $\beta>0$
(and hence $\alpha>0$), we observe that
\begin{equation}\label{sti14b}
  u_n|\Ga = -\frac1{\eta_n}
   \quext{and }\, \bigg\| \frac1{\eta_n} \bigg\|_{L^2(0,T;V\Ga)} \le c,
\end{equation}
the latter bound following from \eqref{sti11}-\eqref{sti12}.
Hence, by standard regularity results for elliptic problems
with Dirichlet boundary conditions
(see, e.g., \cite[Theorem~3.1.5]{BrG}), we obtain
\begin{equation}\label{sti14c}
  \| u_n \|_{L^2(0,T;H^{3/2}(\Omega))}
   \le c.
\end{equation}
Thanks to \eqref{sti11} and \eqref{sti14},
we can apply the trace theorem
\cite[Theorem~2.7.7]{BrG}, which yields
\begin{equation}\label{sti15}
  \| \dn u_n \|_{L^2(0,T;H\Ga)}
   \le c.
\end{equation}
Consequently, using the regularity of the
boundary initial datum \eqref{hp:uteta:reg}
we obtain
\begin{equation}\label{sti16}
  \| \eta_{n,t} \|_{L^2(0,T;H\Ga)}
   + \| \eta_n \|_{L^2(0,T;H^2(\Gamma))}
   \le c.
\end{equation}
Next, we observe that
\begin{equation}\label{sti17}
  \Delta\Ga \Big(-\frac1{\eta_n}\Big)=
   \frac1{\eta_n^2}\Delta\Ga \eta_n
    - 2 \frac{|\nabla\Ga \eta_n|^2}{\eta_n^3}.
\end{equation}
Thus, using the boundary analogue of \eqref{sti12},
\eqref{sti16}, and the Gagliardo-Nirenberg inequality
\begin{equation}\label{GN}
  \| v \|_{W^{1,4}} \le c
   \| v \|_{H^2}^{1/2} \| v \|_{L^\infty}^{1/2}
   + c \| v \|_{L^\infty},
\end{equation}
which gives that $|\nabla\Ga \eta_n|^2\,\in\,L^2(\Gamma)$
uniformly w.r.t. $n$, we readily obtain that
\begin{equation}\label{sti18}
  \bigg\| \frac1{\eta_n} \bigg\|_{L^2(0,T;H^2(\Gamma))} \le c,
\end{equation}
whence \eqref{sti14c} can be improved to
\begin{equation}\label{sti14d}
  \| u_n \|_{L^2(0,T;H^2(\Omega))}
   \le c
\end{equation}
and the same bound holds for $\{\teta_n\}$ thanks again to
\eqref{sti12} and \eqref{GN}.

Let us now deal with the case when $\beta=0$. Then, if
it is also $\alpha = 0$, we are just dealing with
no-flux conditions. Hence, from \eqref{sti14}
we directly deduce \eqref{sti14d}
and, consequently, \eqref{reg:teta:reg3}.
If, instead, $\beta = 0$ and $\alpha >0$, we then have
that $\dn u_n = - \alpha \eta_{n,t}$, so that,
from \eqref{sti13}, \eqref{sti14}, and regularity results
for Neumann elliptic boundary value problems, we deduce
once more \eqref{sti14c}. However, it does not seem possible
to arrive at \eqref{sti14d} since we do not have sufficient
regularity of $\dn u_n$.

In all cases, standard applications of the Aubin-Lions lemma
permit to take the limit $n\nearrow\infty$ and obtain
(up to the extraction of a subsequence) existence of
a triplet $(\teta,\eta,u)$ solving Problem~(P) and complying
with the regularity properties \eqref{reg:teta:reg2}-\eqref{reg:eta:reg2}
and, possibly, \eqref{reg:teta:reg3}.
In particular, it is worth noting that, in the case $\alpha=\beta=0$
there is no boundary function $\eta$. Otherwise, $\eta$ coincides
with the trace of $\teta$ thanks to the quoted regularity
properties and to continuity of the trace operator.
The proof is concluded.
\end{proof}
%

%%%%%%%%%%5%%%%%%%%%%%%%%%%%%%%%%%%%%%%%%%%%%%%%%%%%%%%%%%%%%%%%%%%%%%%%%%%%%%
%%%%%%%%%%5%%%%%%%%%%%%%%%%%%%%%%%%%%%%%%%%%%%%%%%%%%%%%%%%%%%%%%%%%%%%%%%%%%%

\section{Weak solutions}
\label{sec:final}

%%%%%%%%%%5%%%%%%%%%%%%%%%%%%%%%%%%%%%%%%%%%%%%%%%%%%%%%%%%%%%%%%%%%%%%%%%%%%%

\subsection{Existence for finite-energy data}
\label{sec:proof:en}

\bete\label{teo:energy}
 Let assumptions\/ \eqref{hpteta0w}-\eqref{hpeta0w}
 and\/ \eqref{hypfw} hold and let us assume that,
 if $\alpha=0$, then $\beta=0$.
 Then, {\rm Problem~(P)} admits at least one
 energy solution $(\teta,u)$
 satisfying the uniform\/ {\rm energy estimate}
 \begin{equation}\label{en-est}
   \calE(t) + \int_0^t
    \Big( \frac12\|\nabla u\|^2 + \beta \|\nabla\Ga \log \teta\|\Ga^2 \Big)
    \le \EE_0 + c \| f \|_{L^2(0,t;L^{6/5}(\Omega))}^2,
     \quad\perogni t\ge 0.
 \end{equation}
 and the further regularity property
 \begin{equation}\label{LinftyuT}
   \| u(t) \|_{L^\infty(\Omega)}
    \le Q\big(\EE_{0},M_{\epsilon},\tau^{-1}\big)
    \quext{for a.e.~}\,t\ge \tau,~~\tau\in (0,1),
 \end{equation}
 where we have set
 \begin{equation}\label{defi:M}
   M\ee :=
    \| f \|_{L^2(0,+\infty;L^{6/5}(\Omega))}
     + \| f \|_{\calT^2(0,+\infty;L^{3+\epsilon}(\Omega))}.
 \end{equation}
\ente
\begin{proof}
As before, we start with the case when $\alpha$ and
$\beta$ are both strictly positive. The other cases
can be treated with small variants that will be outlined
at the end of the proof. In order to apply
Theorem~\ref{teo:glob:reg}, we consider
a sequence of initial data $\{\teta\zzn\}$ such that,
as $n\nearrow\infty$,
\begin{align} \label{tetazzn}
  & \teta\zzn \in \calV,  \qquad
   0 < \tetagiu_n\le \teta\zzn(x)\le \tetasu_n < +\infty,
    \quext{a.e.~in }\, \Omega~~\text{and on }\,\Gamma,\\
 \label{etazzn}
  & \teta\zzn \to \teta_0 \quext{strongly in }\,L^1(\Omega), \qquad
     \teta\zzn|\Ga \to \eta_0 \quext{strongly in }\,L^1(\Gamma),\\
 \label{etazzn2}
  & \limsup_{n\nearrow\infty} \left( 
     \io \log^- \teta\zzn + \iga \log^- \teta\zzn \right) < \infty. 
\end{align}
The construction of such a
sequence is sketched in Section~\ref{sec:approx} below.
Notice that we do not need to approximate $f$.
Then, for any $n\ge 0$, the assumptions of
Theorem~\ref{teo:glob:reg} are fulfilled and we have
existence of a ``smooth'' solution $(\teta_n,u_n)$.
Moreover, thanks to the energy estimate in 
Lemma~\ref{lemma:energy} (cf.~in particular
\eqref{stimau}), the assumptions of Lemma~\ref{moser-u}
are satisfied over the generic time interval $(t,t+2)$. 
Hence, thanks to \eqref{st-abstr-u}, we obtain 
\begin{equation}\label{energyeq*}
  \| u_n(t) \|_{L^\infty(\Omega)} \le Q(\EE_0,\tau^{-1})
   \quext{for a.e.~}\,t\ge \tau>0,
\end{equation}
where it is worth noting once more that here $Q$ does
not depend on the final time $T$. Indeed, we are
applying Lemma~\ref{moser-u} on a time interval of fixed
length~$2$. Hence, we obtain a bound which is uniform for
large~$T$. 

At this point, we use
$L^1$-techniques in order to take the limit~$n\nearrow\infty$.
To this aim, we work on the generic interval $(0,T)$
and rewrite the approximate equation \eqref{calore}
in the equivalent form (which is possible since
$u_n$ has the good regularity properties
\eqref{reg:u:reg2})
\begin{equation} \label{eq-u}
  \partial_t u_n = u_n^2 \Delta u_n + u_n^2 f.
\end{equation}
Then, we test \eqref{eq-u} with $\fhi\in L^{\infty}(0,T;H^{s}_{0}(\Omega))$,
where we choose $s>\frac{3}{2}$ so
that $H^s_{0}(\Omega)\subset C^0(\barO)$,
continuously. Note that, since $\fhi=0$ on $\Gamma$, we
do not have to take into account the
contribution of the boundary. Using
\eqref{stimau} and \eqref{energyeq*}, it is a standard
matter to get
\begin{equation}\label{tolim-11}
  \| u_{n,t} \|_{L^1(\tau,T;H^{-s})}
   \le c.
\end{equation}
Next, for any $\epsi > 0$ and $r \in \RR$ we introduce
the approximate $\sign$ function as $\sign_\epsi(v):=\frac{v}{\epsi + |v|}$.
Then, we write the approximate equations \eqref{calore}
and \eqref{calorega} for the indexes $m$ and $n$,
take the differences and test them, respectively,
by $\sign_\epsi(u_m- u_n)$ and by
$\sign_\epsi(u_{m,\Gamma} - u_{n,\Gamma})$.
Indeed, for any $m,n\in \NN$ and $\epsi>0$, we have
\begin{equation}\label{tolimi-11}
  \sign_\epsi(u_m - u_n) \in L^2(0,T;V\Ga)
\end{equation}
since $\sign_\epsi$ is Lipschitz.
Then, we arrive at
\begin{align}\no
  & \big( (\teta_{m,t} - \teta_{n,t}), \sign_\epsi(u_m - u_n) \big)
   + \alpha \big( (\eta_{m,t} - \eta_{n,t}), \sign_\epsi(u_m - u_n) \big)\Ga\\
  \label{giulio-11}
  & \mbox{}~~~~~
   - \beta \big( \Delta\Ga (\eta_m - \eta_n), \sign_\epsi(u_m - u_n) \big)\Ga
    \le 0.
    %\le \big( f_n - f_m, \sign_\epsi(u_m - u_n) \big)
    %\le \| f_n - f_m \|_1.
\end{align}
%
%where we intend, at this level, that $\alpha=\beta=1$.
Here, we used integration by parts to prove that
\begin{equation}\label{tolimi-12}
  - \big( \Delta (u_m - u_n), \sign_\epsi(u_m - u_n) \big)
   \ge - \big( \dn (u_m - u_n), \sign_\epsi(u_m - u_n) \big)\Ga,
\end{equation}
which is possible thanks to the good regularity of approximate
solutions, to monotonicity of $\sign_\epsi$,
and to \eqref{tolimi-11}. Thus, integrating
\eqref{tolimi-11} over $(0,t)$, for a generic $t\in(0,T]$,
we obtain
\begin{align}\no
  & \itt \big( (\teta_{m,t} - \teta_{n,t}), \sign_\epsi(u_m - u_n) \big)
   + \alpha \itt \big( (\eta_{m,t} - \eta_{n,t}), \sign_\epsi(u_m - u_n) \big)\Ga\\
  \label{giulio-12}
  & \mbox{}~~~~~
   - \beta \itt \big( \Delta\Ga (\eta_m - \eta_n), \sign_\epsi(u_m - u_n) \big)\Ga
    \le 0.
  %  \le \| f_n - f_m \|_{L^1(0,T;L^1(\Omega))}.
\end{align}
Next, we take the limit $\epsi\searrow 0$ and notice that,
since $v\mapsto -1/v$ is monotone increasing,
it turns out that $\sign_0(\teta_m -\teta_n)=\sign_0(u_m -u_n)$
almost everywhere in $(0,T)\times \Omega$
and $\sign_0(\eta_m - \eta_n)=\sign_0(u_m -u_n)$
almost everywhere in $(0,T)\times \Gamma$,
where $\sign_0$ denotes the sign-like function
with $\sign_0(0)=0$. Then,
\begin{align}\no
  & \itt \big( (\teta_{m,t} - \teta_{n,t}), \sign_0(\teta_m - \teta_n) \big)
   + \alpha \itt \big( (\eta_{m,t} - \eta_{n,t}), \sign_0(\eta_m - \eta_n) \big)\Ga\\
  \label{giulio-13}
  & \mbox{}~~~~~
   - \beta \itt \big( \Delta\Ga (\eta_m - \eta_n), \sign_0(\eta_m - \eta_n) \big)\Ga
    \le 0.
   % \le \| f_n - f_m \|_{L^1(0,T;L^1(\Omega))}.
\end{align}
Applying the Brezis-Strauss lemma (see \cite[Lemma 2]{BreStra}),
we can integrate by parts the boundary Laplacian,
which gives a nonnegative contribution. Thus, we arrive at
\begin{equation}\label{giulio-14}
  \big\| \teta_m(t) - \teta_n(t) \big\|_1
   + \alpha \big\| \eta_m(t) - \eta_n(t) \big\|_{1,\Gamma}
   \le \big\| \teta\zzm - \teta\zzn \big\|_1
   + \alpha \big\| \teta\zzm - \teta\zzn \big\|_{1,\Gamma}.
\end{equation}
%
%\begin{align}\no
%  & \big\| \teta_m(t) - \teta_n(t) \big\|_1
%   + \alpha \big\| \eta_m(t) - \eta_n(t) \big\|_{1,\Gamma}\\
% \label{giulio-14}
%  & \mbox{}~~~~~
%   \le \big\| \teta\zzm - \teta\zzn \big\|_1
%   + \alpha \big\| \teta\zzm - \teta\zzn \big\|_{1,\Gamma}.
%   + \| f_n - f_m \|_{L^1(0,T;L^1(\Omega))}.
%\end{align}
%
Taking the supremum for $t\in[0,T]$, we finally obtain
\begin{equation}\label{tolim-12}
  \teta_n \to \teta \quext{strongly in }\,C^0([0,T];L^1(\Omega)), \qquad
    \alpha \eta_n \to \alpha \eta \quext{strongly in }\,C^0([0,T];L^1(\Gamma)),
\end{equation}
for suitable limit functions $\teta$ and $\eta$. Moreover, using
\eqref{stimau}, \eqref{tolim-11}, and the Aubin-Lions lemma,
we obtain that, up to the extraction of a subsequence,
\begin{equation} \label{convergence1}
  u_n \to u \quext{weakly in }\,L^2(0,T;V)
  \qquext{and strongly in }\,L^2(\tau,T;H^{1-\delta}(\Omega)),
\end{equation}
for all $\tau\in(0,1)$, $\delta\in(0,1)$,
and for some limit function $u$.
In particular, taking a test function $\xi$
satisfying \eqref{regotest}, we can test both
\eqref{calore} and \eqref{calorega} (at the $n$-level)
by $\xi$ and integrate by parts. This gives
the $n$-analogue of \eqref{caloreveryw}.
Moreover, the convergence properties proved
above permit to take the limit $n\nearrow\infty$
to show that \eqref{caloreveryw} holds
also at the limit level. Hence, to conclude
the proof it remains to identify the functions
$\teta$ and $\eta$ in terms of $u$. To do this,
we first observe that, combining the first
\eqref{tolim-12} with \eqref{convergence1},
it follows that $\teta=-1/u$
almost everywhere in $(0,T)\times \Omega$.

Next, by \eqref{convergence1} and
continuity of the trace operator, we have,
for all $\tau\in(0,T)$,
\begin{equation}\label{convergence1gamma}
  (u_n)_{\Gamma} \to u_\Gamma
   \quext{strongly in } \,L^2(\tau,T;H^{1/2-\delta}(\Gamma)),
\end{equation}
with $u_\Gamma$ denoting the trace of $u$. In particular,
thanks to arbitrariness of $\tau$,
pointwise convergence holds on $(0,T)\times\Gamma$,
up to the extraction of a further subsequence.
Being $\eta_n = -1/(u_n)\Ga$
almost everywhere in $(0,T)\times \Gamma$ and for
all $n\in\NN$, we then deduce that
\begin{equation}\label{convergence2gamma}
  \eta_n \to - \frac1{u\Ga}
   \quext{a.e.~on } (0,T)\times \Gamma.
\end{equation}
Combining this with the second \eqref{tolim-12},
we then obtain relation \eqref{trace}. This concludes the
proof of the Theorem in the case $\alpha>0$, $\beta >0$.

\medskip

Let us now consider the case $\alpha > 0$ and $\beta = 0$
(the case $\alpha=\beta=0$ is essentially already known
and, in any case, it is simpler to treat).
Since Lemma \ref{moser-u} does not depend on
whether $\beta$ is zero or not, estimate \eqref{energyeq*}
(and consequently \eqref{tolim-11}) still holds when
$\beta=0$. Moreover, it is easy to check that the Cauchy
argument \eqref{tolimi-11}-\eqref{tolim-12}
can be reproduced also when $\beta=0$.
The proof of the Theorem is concluded.
\end{proof}
\beos\label{boundarylayer}
 Note that, in the case when $\alpha>0$ and $\beta>0$,
 $\eta_0$ need not be the trace of $\teta_0$.
 Thus, there is a boundary layer in the sense specified
 by \eqref{convergence1gamma}-\eqref{convergence2gamma}.
 However, \eqref{trace} still makes sense
 thanks to instantaneous regularization properties.
\eddos

%%%%%%%%%%5%%%%%%%%%%%%%%%%%%%%%%%%%%%%%%%%%%%%%%%%%%%%%%%%%%%%%%%%%%%%%%%%%%%

\subsection{Approximation of data}
\label{sec:approx}

In this part, we sketch the approximation of data
needed in the proofs of Theorems~\ref{teo:glob:reg}
and \ref{teo:energy}. Since the procedures are rather
standard, we just give the highlights without entering
too much into details.

\smallskip

\noindent%
{\bf Approximation of data for Theorem \ref{teo:glob:reg}.}~~%
We start with the initial data.
Let us given $\teta_0 \in \calV$ such that
$\tetagiu\le \teta_0(x)\le \tetasu$ for all $x\in\barO$
(as in the statement, see \eqref{hp:uteta:reg}),
and let $\eta_0$ be its trace on $\Gamma$.
Since we want to apply Theorem~\ref{teo:loc:reg}
(cf., in particular, \eqref{hpteta:reg}),
for any $n\in\mathbb{N}$ we need to have an approximate
datum $\teta_{0,n}$ (with trace $\eta_{0,n}$)
such that $\teta_{0,n}\in H^2(\Omega)$ for all $n$
and
\begin{enumerate}
 \item $\teta_{0,n}\xrightarrow{n \nearrow +\infty}\teta_{0}$ in $\mathcal{V}$,
 \item $\tetagiu\le \teta_{0,n}(x)\le \tetasu$ for all $x\in\barO$.
\end{enumerate}
To construct $\teta_{0,n}$ we then consider the
elliptic problem
\begin{equation} \label{approxdataTh1}
 \displaystyle\begin{cases}
  \teta_{0,n}-\frac{1}{n}\Delta \teta_{0,n}= \teta_0~~\text{in }\,\Omega,\\
   -\frac{1}{n}\partial_n\teta_{0,n}
       = -\frac{1}{n}\Delta\Ga\eta_{0,n}
         +\eta_{0,n}-\eta_{0}~~\text{on }\,\Gamma,~~
      \teta_{0,n}|_\Gamma = \eta_{0,n}.
 \end{cases}
\end{equation}
Then, the standard elliptic theory gives that, for any
fixed $n$, $\teta_{0,n}$ lies in $H^2(\Omega)$
and $\eta_{0,n}$ lies in $H^2(\Gamma)$.

To prove the convergence $\teta_{0,n}\to \teta_0$,
we test \eqref{approxdataTh1}$|_{\Omega}$
by $\big(\teta_{0,n} -\Delta \teta_{0,n}\big)$
and, correspondingly, \eqref{approxdataTh1}$|_{\Gamma}$
by $(1+n^{-1})\eta_{0,n}-\eta_0$.
Then, by integrations by parts and Young's inequality,
it is not difficult to obtain
\begin{equation*}
  \|\teta_{0,n}\|^{2}_{\mathcal{V}}
   + 2n \|\eta_{0,n}-\eta_{0}\|^{2}_{\Gamma}
  \le \|\teta_{0}\|^{2}_{\mathcal{V}},
\end{equation*}
which clearly implies the desired convergence. On the other hand,
to prove uniform boundedness of $\teta_{0,n}$,
we use a maximum principle argument.
We prove only the upper bound, the lower one being completely
equivalent. For a generic $M>\tetasu$,
we test $\eqref{approxdataTh1}|_{\Omega}$
by $(\teta_{0,n}-M)^{+}$. Integrating by parts the
Laplacian and using $\eqref{approxdataTh1}|_{\Gamma}$,
we then obtain
\begin{eqnarray}\no
  \int_{\Omega}\teta\zzn(\teta\zzn -M)^+
   + \int_{\Gamma}\eta\znn(\eta\znn-M)^+ \le \int_{\Omega}\teta_0(\teta\zzn -M)^+
   + \int_{\Gamma}\eta_0(\eta\znn-M)^+ .
\end{eqnarray}
Adding to both sides the quantity
$-M\int_{\Omega}(\teta\zzn-M)^+ -M\int_{\Gamma}(\eta\zzn-M)^+$,
we arrive at
\begin{equation}\label{antonio11}
  \|(\teta\zzn -M)^+\|^2
   + \|(\eta\zzn -M)^+\|^2_{\Gamma}
  \le \int_{\Omega} (\teta_0-M)(\teta\zzn -M)^+
   + \int_{\Gamma}(\eta_0-M)(\eta\znn-M)^+
   \le 0.
\end{equation}
This clearly gives the desired inequality $\teta\znn(x) \le\tetasu$
for all $x\in\barO$.

\smallskip

Concerning the source term, recalling that $f$ satisfies \eqref{hypfw},
a combination of truncation and mollification techniques,
together with a suitable correction of the spatial mean values
permits to construct a family of functions $\{f_n\}$, $n\in \NN$,
such that
\begin{align} \label{deffn}
  & f_n \in C^0([0,T]\times \barO), \qquad
   \| f_n \|_{L^\infty((0,T)\times \Omega)} \le n,\\
 \label{deffn2*}
  & f_n \to f \quext{strongly in }\,L^1((0,T)\times \Omega),\\
 \label{deffn3*}
  & (f_n(t))\OO = 0 \quext{for all }\,t\in[0,T],\\
 \label{deffn4*}
  & \| f_n \|_{L^2((0,T;L^{3+\epsilon}(\Omega))} \le c.
\end{align}
The details of this construction are left to the reader.
Notice that \eqref{deffn}-\eqref{deffn4*} suffice
both to apply Theorem~\ref{teo:loc:reg} at the level $n$
to get a local smooth solution and to perform the
estimates of Theorem~\ref{teo:glob:reg} uniformly in
$n$ in order to let $n\to \infty$.
%
% Moreover,
%we approximate the maximal monotone graph $-1/(\cdot)$,
%by taking an approximating family $\{\gamma_n\}$
%where $\gamma_n$ are strictly monotone functions
%such that
%%
%\begin{align} \label{gamman}
%  & \gamma_n\in C^2(\RR;\RR),~~\qquad
%   \gamma_n,~\gamma_n',~\gamma_n^{-1}\in L^\infty(\RR),\\
%  \label{gamman2}
%  & \gamma_n(r) = -\frac1{r} \quext{for all }\,
%    r\in \Big[\frac1{2n},2n\Big].
%\end{align}
%
%Also the construction of this family is standard. For instance, one
%can take $\gamma\,\equiv\,\gamma_n$ for $r\,\in\,\Big[\frac1{2n},2n\Big]$
%and $\gamma_n$ equal to the Taylor expansion of order one
%outside $\Big[\frac1{2n},2n\Big]$. Notice
%that, of course, the Lipschitz constants of $\gamma_n$
%and $\gamma_n^{-1}$ will explode as $n\nearrow\infty$.

\smallskip

\noindent%
{\bf Approximation of data for Theorem \ref{teo:energy}.}~~%
We detail such an approximation just in the case $\alpha>0$,
which is a little bit trickier. So, let
$(\teta_0,\eta_0) \in L^1(\Omega)\times L^1(\Gamma)$
such that $\log \teta_0 \in L^1(\Omega)$ and $\log \eta_0 \in L^1(\Gamma)$.
We need to construct $\teta_{0,n}$ in such a way that
properties \eqref{tetazzn}-\eqref{etazzn} hold.
Then, we consider first the function $\teta_0$. Using that $\Gamma$ is
smooth, we can first extend it (e.g., by reflection) to a neighbourhood
$\Omega_*$ of $\barO$. It is then clear that the new function,
note it as $\teta_{0,*}$, lies in $L^1(\Omega_*)$; moreover,
$\log \teta_{0,*} \in L^1(\Omega_*)$. Next, we truncate
$\teta_{0,*}$, setting
\begin{equation} \label{approx-initial-cond}
  \teta_{0,n}^{(1)}:=\min\left\{\max\left\{\teta_{0,*},\frac{1}{n}\right\}, n\right\},
   \quext{a.e.~in }\, \Omega_*.
\end{equation}
Finally, we regularize, setting
$\teta_{0,n}^{(2)}:= \rho_n * \teta_{0,n}^{(1)}$, where
$\{\rho_n\}$ is a suitable sequence of smooth
and compactly supported mollifiers.
Then, straightforward checks (based on the properties
of convolutions and on Lebesgue's theorem)
and a standard diagonal argument
permit to verify that $\teta_{0,n}^{(2)}$ is smooth and
tends to $\teta_0$ strongly in $L^1(\Omega)$. Moreover,
based on Jensen's inequality, it is not difficult
to verify that, for some $c>0$ independent of $n$,
\begin{equation} \label{log-approx}
  \io \log^- \teta_{0,n}^{(2)}
   \le c \Big(1 + \io \log^- \teta_0 \Big).
\end{equation}
Finally, we pass to the boundary component. First of all,
thanks to smoothness of $\Gamma$, we can find $\epsilon>0$
such that $\eta_0$ can be extended (e.g., constantly
along directions orthogonal to $\Gamma$) to a function
$\eta_{0,*}$ defined on a neighbourhood
$\Gamma_\delta:=\{x\in\RR^3 : d(x,\Gamma)\le \delta \}$
of $\Gamma$. Thanks to Fubini's theorem, it is then clear that
both $\eta_{0,*}$ and $\log\eta_{0,*}$
lie in $L^1 (\Gamma_\delta)$. Then, we truncate $\eta_{0,*}$
(as in \eqref{approx-initial-cond})
obtaining $\eta_{0,n}^{(1)}$ (which can be seen as a function
defined on the whole of $\RR^3$). Next, we
mollify $\eta_{0,n}^{(1)}$, introducing
$\eta_{0,n}^{(2)}:= \rho_n * \eta_{0,n}^{(1)}$,
for $\rho_n$ supported, say, on the ball
$\overline{B}(0,1/n)$.
Finally, we take a cutoff function
$\psi_n\in C^\infty(\RR^3;[0,1])$ such
that $\psi_n$ is identically one on $\Gamma_{1/2n}$ and
$\psi_n$ is supported on $\Gamma_{1/n}$. Then,
we set $\eta_{0,n}^{(3)}:=\eta_{0,n}^{(2)}\psi_n$
in such a way that $\eta_{0,n}^{(3)}$ belongs to
$C^\infty(\RR^3)$ and $\eta_{0,n}^{(3)}$ tends
to $0$ in $L^1(\Omega)$, while its trace
tends to $\eta_0$ in $L^1(\Gamma)$. Moreover,
as above, one can check that
$\log^-(\eta_{0,n}^{(3)})$ is uniformly
controlled in $L^1(\Gamma)$ in the sense of
\eqref{log-approx}. Then, the required
approximation of $\teta_0$ is obtained simply taking 
$\teta_{0,n}:=(1-\psi_n) \teta_{0,n}^{(2)}+\eta_{0,n}^{(3)}$. 
%
%
%Finally, to regularize $f$, we can simply truncate it
%at the levels $\pm n$ to obtain
%a bounded $f_n$ such that
%
%\begin{equation} \label{fL1}
%  f_n \to f \quext{strongly in }\,
%   L^1((0,T)\times\Omega).
%\end{equation}
%

%%%%%%%%%%%%%%%%%%%%%%%%%%%%%%%%%%%%%%%%%%%%%%%%%%%%%%%%%%%%%%%%%%%%%%%%%%%%%%%%%%%%%%%%%%%%%%%%%%%%%%
%%%%%%%%%%%%%%%%%%%%%%%%%%%%%%%%%%%%%%%%%%%%%%%%%%%%%%%%%%%%%%%%%%%%%%%%%%%%%%%%%%%%%%%%%%%%%%%%%%%%%%

\subsection{Solutions with regularizing effects for $\teta$}
\label{subsec:asymteta}

In this last part, we extend the previous results in 
three directions.
First, we prove that, if $\teta_0$ (and, possibly, $\eta_0$)
enjoy higher summability properties, then
there exist weak solutions whose component $\teta$
satisfies time-regularization properties in the spirit of
Lemmas~\ref{moser-teta-2} and~\ref{moser-teta-2b}.
Second, we demonstrate that, under the same type of
conditions on the initial data, existence holds also for
$\alpha=0$ and $\beta>0$ (recall that we could not deal
with this case for $L^1$ initial data, 
cf.~Theorem~\ref{teo:energy}).
Third, we see that uniqueness holds in the class
of solutions with regularizing effects. We start
with analyzing regularity:  
\bete\label{teo:smoothing}
 Let assumptions\/ \eqref{hpteta0w}-\eqref{hpeta0w}
 and \eqref{hypfw} hold
 and let in addition, for some $\epsilon\in(0,1)$,
 \begin{equation}
 \label{hypf2}
   f \in L^2(0,+\infty;L^{3+\epsilon}(\Omega)), \qquad
    N\ee:= \| f \|_{L^2(0,+\infty;L^{3+\epsilon}(\Omega))}.
 \end{equation}
 Moreover, if $\alpha>0$ and $\beta=0$, let
 \begin{equation}\label{hpteta0-s}
   \teta_0 \in L^{4+\epsilon}(\Omega), \qquad
    \alpha \teta_0 \in L^{4+\epsilon}(\Gamma),
 \end{equation}
 whereas in the other cases let
 \begin{equation}\label{hpteta0}
   \teta_0 \in L^{3+\epsilon}(\Omega), \qquad
    \alpha\teta_0 \in L^{3+\epsilon}(\Gamma).
 \end{equation}
 Then, {\rm Problem~(P)} admits at least one
 energy solution $(\teta,u)$
 satisfying~\eqref{en-est}, \eqref{LinftyuT},
 together with the regularization estimate
 \begin{equation}\label{st:uniform}
   \| \teta(t) \|_{L^\infty(\Omega)}
    \le Q\big(\EE_{0+\epsilon},N\ee,\tau^{-1}\big)
    \quad \perogni t\ge\tau,~~\tau\in(0,1),
 \end{equation}
 where we have set
 \begin{equation}\label{defi:E0ee}
   \mathbb{E}_{0+\epsilon}
   :=  \calE(\teta_0)
    + \|\teta_0\|_{L^{3+\epsilon}(\Omega)}^{3+\epsilon}
    + \alpha\|\teta_0 \|_{L^{3+\epsilon}(\Gamma)}^{3+\epsilon},
 \end{equation}
 the exponents $3+\epsilon$ being all replaced by $4+\epsilon$
 in the case when $\alpha>0$ and $\beta=0$. Moreover,
 in all cases with the exception of $\alpha=0$ and $\beta>0$,
 we have
 \begin{align}\no
   & \| \teta_t \|_{L^2(t,t+1;H)}
    + \alpha \| \eta_t \|_{L^2(t,t+1;H\Ga)}
    + \| u \|_{L^\infty(t,+\infty;V)}
    + \beta \| \eta \|_{L^\infty(t,\infty;V\Ga)}\\
  \label{st:uniform2}
   & \mbox{}~~~~~~~~~~
   \le Q\big(\EE_{0+\epsilon},N_{\epsilon},\tau^{-1})
    \quad\perogni t\ge\tau,~~\tau\in (0,1).
 \end{align}
 Moreover, estimates\/ \eqref{st:uniform} and\/ \eqref{st:uniform2} 
 are uniform for large values of $t$. In other words, the expression
 of $Q$ is independent of $t$.  
\ente
\begin{proof}
For the sake of simplicity, we just prove the theorem
by directly working on the ``limit'' solutions without
referring to an explicit approximation scheme.
That said, we first observe that the energy estimate
\eqref{stimau} still holds. Moreover, we can still
rely on the conclusion of Lemma~\ref{moser-u}.
Next, thanks to assumption \eqref{hypf2},
we can apply Lemma~\ref{moser-teta-1} over the generic
time interval $(0,T)$. Note that the estimates provided
by Lemma~\ref{moser-teta-1} are uniform with respect to~$T$.
We then obtain 
%
%({\bf it should be enough to assume $L^3$-valued summability
%and $L^{3+\epsilon}$-local summability, but I am
%not sure it is worth adding technicalities into
%the proof in order to take care of such a small extension})
%
\begin{equation} \label{fine-11}
  \| \teta(t) \|_{L^{3+\epsilon}(\Omega)}
   + \alpha \| \teta(t) \|_{L^{3+\epsilon}(\Gamma)}
  \le Q\big( N\ee, \EE_{0+\epsilon}\big)
\end{equation}
(here and below, $3+\epsilon$ is replaced by $4+\epsilon$
in case $\alpha>0$ and $\beta=0$). Hence, for $t\ge \tau > 0$, we can
apply Lemma~\ref{moser-teta-2} (or Lemma~\ref{moser-teta-2b})
over the generic time interval $(t,t+2)$, which has
fixed finite length. This gives
\eqref{st:uniform}.

Finally, in all cases with the exception of $\alpha=0$ and $\beta>0$,
we can apply Lemma~\ref{lemma-time} over the generic
time interval $(t,t+1)$ where $t\ge \tau>0$, which gives
\eqref{st:uniform2}.

Then, to conclude the proof, it just remains to show
that, in the case when $\alpha=0$ and $\beta>0$
(that we set equal to $1$ for simplicity),
a weak solution still exists under the above
assumptions. To this aim, we consider the system
\begin{align}\label{eq-abstrn}
  & \teta_{n,t} - \Delta u_n = f, \quad \teta_n = -1/u_n,
   \quext{in }\,\Omega,\\
 \label{bc-abstrn}
  & \frac1n \eta_{n,t} - \Delta\Ga \eta_n = - \dn u,
  \quext{on }\,\Gamma,
\end{align}
complemented with the usual initial conditions.
Then, for all $n\in\NN$, there exists at least one weak
solution $(\teta_n,u_n)$.
Moreover, thanks to \eqref{st:uniform2} and
to regularity arguments similar to those performed in
Section~\ref{sec:global}, $(\teta_n,u_n)$ is smooth enough
in order for the system to make sense in the above
``strong'' form, at least on time intervals of
the form $(\tau,T)$ for all $\tau>0$.

In addition to that, we still have
estimates \eqref{en-est} and \eqref{LinftyuT}. Moreover,
it is worth noting that \eqref{LinftyuT} holds
{\sl independently of $n$}. Indeed, looking back
at the proof of Lemma~\ref{moser-u}, it is immediate
to check that estimate \eqref{st-abstr-u} is
independent of $\alpha$. Actually,
when one performs the iteration argument
(see \eqref{moser11}) one simply has a functional
$J_i$ that depends on $\alpha$, but the bound for such
a functional remains unchanged. The same holds when
applying Lemma~\ref{moser-teta-2}. In conclusion,
we have the uniform bound
\begin{equation}\label{sti:31}
  \| \teta_n(t) \|_{L^\infty(\Omega)}
   + \| u_n(t) \|_{L^\infty(\Omega)}
   \le Q\big(\EE_{0+\epsilon},N\ee,\tau^{-1}\big)
   \quad \perogni t\ge\tau,~~\tau\in(0,1).
\end{equation}
Thus, thanks also to Lemma~\ref{moser-teta-1},
for any $\tau\in(0,1)$, $T\ge\tau$, we have
\begin{equation} \label{fine-21}
  \teta_n \to \teta \quext{weakly star in }\,L^\infty(0,T;L^{3+\epsilon}(\Omega))
   \quext{and weakly star in }\,L^\infty((\tau,T)\times\Omega).
\end{equation}
Moreover, since $\beta>0$, as an additional consequence of estimate
\eqref{lemma2-11} (with $p=3+\epsilon$) we have
\begin{equation} \label{fine-22}
  \big\| \teta_n^{\frac{1+\epsilon}2} \big\|_{L^2(0,T;V)}
  + \big\| \eta_n^{\frac{3+\epsilon}2} \big\|_{L^2(0,T;V\Ga)}
  \le c.
\end{equation}
In particular, being
\begin{equation} \label{fine-23}
  \nabla \teta_n = \frac{2}{1+\epsilon} \teta_n^{\frac{1-\epsilon}2}
         \nabla \teta_n^{\frac{1+\epsilon}2},
\end{equation}
we have, from \eqref{fine-21}-\eqref{fine-22},
\begin{equation} \label{fine-24}
  \| \nabla \teta_n \|_{L^2(0,T;L^{\frac{3+\epsilon}{2}}(\Omega))}
   \le c \big \| \nabla \teta_n^{\frac{1+\epsilon}2} \big\|_{L^2(0,T;H)}
    \| \teta_n^{\frac{1-\epsilon}2} \|_{L^\infty(0,T;L^{\frac{6+2\epsilon}{1-\epsilon}}(\Omega))}
   \le c.
\end{equation}
In particular, this fact tells us that,
in the present regularity setting, $\eta_n$ can be directly
seen as the trace of $\teta_n$. More precisely,
applying the trace theorem, we have
\begin{equation} \label{fine-24b}
  \| \eta_n \|_{L^2(0,T;W^{\frac{1+\epsilon}{3+\epsilon},\frac{3+\epsilon}{2}}(\Gamma))}
   \le c.
\end{equation}
On the other hand, testing \eqref{eq-abstr} by a generic
function $\phi\in H^1_0(\Omega)$ of unit norm and recalling
\eqref{en-est} and \eqref{hypf2}, we obtain
\begin{equation} \label{fine-25}
  \teta_{n,t} \to \teta_t \quext{weakly in }\,L^2(0,T;H^{-1}(\Omega)).
\end{equation}
Thus, using \eqref{fine-24}, \eqref{fine-25} and
the Aubin-Lions lemma, we infer
\begin{equation} \label{fine-26}
  \teta_n \to \teta \quext{strongly in }\,L^2(0,T;W^{1-\sigma,\frac{3+\epsilon}{2}}(\Omega))
\end{equation}
for all $\sigma\in(0,1)$. Hence, by continuity of the trace operator,
\begin{equation} \label{fine-26a}
  \eta_n \to \eta=\teta|\Ga \quext{strongly in }\,L^2(0,T;L^p(\Gamma))
   \quext{for some }\,p>1.
\end{equation}
Thus, writing the weak formulation \eqref{caloreveryw} of system
\eqref{eq-abstrn}-\eqref{bc-abstrn}, it is immediate to
check that one can take the limit $n\to \infty$ therein.
In particular, it can be standardly proved that
\begin{equation} \label{fine-26b}
  u_n \to u = -1/\teta \quext{strongly in }\,L^p(0,T;L^p(\Omega)),
   \quext{say, for all }\, p\in [1,2).
\end{equation}
Hence, in place of \eqref{trace}, we can now directly
write $\eta=\teta|\Ga$ in the sense of traces and a.e.~in~$(0,T)$.
We may also notice that, thanks to the additional regularity properties
coming from \eqref{fine-24}-\eqref{fine-24b}, we could even
relax a bit the requirements \eqref{regotest}
on the test function $\xi$ (we omit the details).
The proof is concluded.
\end{proof}
Finally, let us come to uniqueness:
\bete\label{teo:newuniq}
 Let\/ \eqref{hpteta0w}-\eqref{hpeta0w} and \eqref{hypfw} hold.
 Let $(\teta_1,u_1)$ and $(\teta_2,u_2)$ be a couple of weak
 solutions to\/ {\rm Problem (P)} over some interval $(0,T)$
 in the sense of\/
 {\rm Definition~\ref{def:weak}}, both satisfying
 \eqref{st:uniform} and \eqref{st:uniform2} and 
 emanating from the same initial datum $\teta_0$. 
 Then $(\teta_1,u_1)$ and $(\teta_2,u_2)$ coincide
 over $(0,T)$.
\ente
\begin{proof}
In view of the fact that $(\teta_1,u_1)$ and $(\teta_2,u_2)$ 
are smooth for strictly positive
times, we can proceed as in Section~\ref{sec:proof:en}
testing the difference of the equations
by $\sign_\epsi(u_1-u_2)$. We then integrate over
$(\tau,T)$ for $\tau>0$ and arrive at the
analogue of \eqref{giulio-12}, namely
\begin{align}\no
  & \int_\tau^t \big( (\teta_{1,t} - \teta_{2,t}), \sign_\epsi(u_1 - u_2) \big)
   + \alpha \int_\tau^t \big( (\eta_{1,t} - \eta_{2,t}), \sign_\epsi(u_1 - u_2) \big)\Ga\\
  \label{giulio-12x}
  & \mbox{}~~~~~
   - \beta \int_\tau^t \big( \Delta\Ga (\eta_1 - \eta_2), \sign_\epsi(u_1 - u_2) \big)\Ga
   \le 0.
\end{align}
Note that we cannot integrate directly over $(0,t)$ since
\eqref{st:uniform} and \eqref{st:uniform2} do not extend
to $\tau=0$; in other words, we do not have sufficient regularity
to use $\sign_\epsi(u_1 - u_2)$ as a test function
over $(0,T)$. Hence, we first need to
take the limit $\epsi\searrow0$, obtaining
\begin{equation}\label{giulio-14x}
  \big\| \teta_1(t) - \teta_2(t) \big\|_1
   + \alpha \big\| \eta_1(t) - \eta_2(t) \big\|_{1,\Gamma}
   \le \big\| \teta_1(\tau) - \teta_2(\tau) \big\|_1
   + \alpha \big\| \eta_1(\tau) - \eta_2(\tau) \big\|_{1,\Gamma}.
\end{equation}
Letting $\tau\searrow 0$ and noting that energy solutions
are continuous with values in $L^1$ (cf.~\eqref{regotetaw}),
we obtain the assertion. The proof is concluded.
\end{proof}
\beos\label{on:uniq}
 In principle, the above proof does not rely directly on hypothesis
 \eqref{hpteta0} (or \eqref{hpteta0-s}). In other words, 
 the conditions on initial data 
 assumed in Theorem~\ref{teo:energy} may suffice.
 However, in that case the uniqueness statement may be vacuous since
 we do not know whether there exist energy solutions satisfying
 \eqref{st:uniform} and \eqref{st:uniform2}. Actually, even
 when \eqref{hpteta0} (or \eqref{hpteta0-s}) holds,
 the proved properties does not exclude that there might
 exist other energy solutions of Problem~(P)
 that {\sl do not}\/ regularize with respect to time.
 The same observation can also be referred to the case
 when $\alpha=0$ and $\beta>0$ since \eqref{st:uniform2}
 is not known to hold under these conditions.
\eddos
%

%%%%%%%%%%%%%%%%%%%%%%%%%%%%%%%%%%%%%%%%%%%%%%%%%%%%%%%%%%%%%%%%%%%%%%%%%%%%%%%%%%%%%%%%%%%%%%%%%%%%%%
%%%%%%%%%%%%%%%%%%%%%%%%%%%%%%%%%%%%%%%%%%%%%%%%%%%%%%%%%%%%%%%%%%%%%%%%%%%%%%%%%%%%%%%%%%%%%%%%%%%%%%

%\section{References}

%\vspace{1cm}

\end{document}